\documentclass[11pt,a4paper,leqno,oneside]{article}

\usepackage{amsmath}
\usepackage{amsfonts}
\usepackage{amssymb}
\usepackage{amsthm}
\usepackage{makeidx}
\usepackage{graphicx}
\usepackage{pb-diagram}
\usepackage{epstopdf}
\usepackage{fancyhdr}
\usepackage{enumitem}
\usepackage[all]{xy}

\newtheorem{teo}{Theorem}%definizioni che danno luogo ai nomi di teoremi, proposizioni, ecc.
\newtheorem{prop}[teo]{Proposition} %la quadra [teo] indica che la numerazione  subordinata a quella dei teoremi
\newtheorem{lem}[teo]{Lemma}
\newtheorem{cor}[teo]{Corollary}
\theoremstyle{definition}
\newtheorem{defin}[teo]{Definition}

\newtheorem{nt}[teo]{Notation}

\newtheorem{nota}[teo]{Note}

\newtheorem{con}[teo]{Conjecture}
\newtheorem{prob}[teo]{Problem}

\newtheorem{rem}[teo]{Remark}

   {\begin{verse}%
     \small }%
   {\end{verse}}

\providecommand{\A}{\mathcal{A}} 
 \providecommand{\U}{\mathcal{U}}

\title{A SURVEY ON CONNES' EMBEDDING CONJECTURE}

\author{VALERIO CAPRARO}

\begin{document}

\maketitle

\setlength{\parskip}{1ex plus 0.5ex minus 0.2ex}

\tableofcontents

\section{Introduction}
In his famous paper \cite{Co}, A. Connes formulated a conjecture
which is now one of the most important open problem in Operator
Algebras. This importance comes from the works of many
mathematicians (above all Kirchberg \cite{Ki}, but also Brown
\cite{Br}, Collins-Dykema \cite{Co-Dy}, Haagerup-Winslow
\cite{Ha-Wi1},\cite{Ha-Wi2}, Radulescu \cite{Ra1},\cite{Ra2},
Voiculescu \cite{Vo2} and many others) who have found some
unexpected equivalent statements showing as this conjecture is
transversal to
almost all the sub-specialization of Operator Algebras.\\
In this survey I would like to give a more or less detailed description of
all these approaches.\\
In the second chapter I am going to recall, more or less briefly,
some preliminary notions (ultrafilters, ultraproducts ...) and give
the original formulation of the conjecture. In the third one I am
going to describe Radulescu's algebraic approach via hyperlinear
groups. In the forth one I am going to describe Haagerup-Winslow's
topological approach via Effros-Marechal topology. In the fifth one
I am going to describe Brown's theorem which connects Connes'
embedding conjecture with Lance's weak expectation property. In the
sixth one I am going to describe briefly other approaches.
\section{Preliminary notions and
original formulation of the conjecture} The most important
preliminary notions are those of ultrafilter and ultraproduct. We
recall the following
\begin{defin}\label{ultrafilter}
Let $X$ be a set and $\U$ a non-empty family of subsets of $X$. We
say that $\mathcal U$ is an ultrafilter if the following properties
are satisfied:
\begin{enumerate}
\item $\emptyset\notin\mathcal U$
\item $A,B\in\mathcal U$ implies $A\cap B\in\mathcal U$
\item $A\in\U$ implies $B\in\U,\forall B\supseteq A$
\item For each $A\subseteq X$ one has either $A\in\U$ or $X\setminus
A\in\U$
\end{enumerate}
\end{defin}
Ultrafilters are very useful in topology, since they can be thought
as a dual notion of \emph{net}, allowing to speak about convergence
in a very general setting. In this overview we are interested only
in the concept of limit along an ultrafilter of a family of real
numbers.
\begin{defin}
Let $\{x_a\}_{a\in A}$ be a family of real numbers and $\U$ an
ultrafilter on $A$. We say that $lim_{\U}x_a=x\in\mathbb R$ if for
any $\varepsilon>0$ one has
$$
\{a\in A:|x_a-x|<\varepsilon\}\in\U
$$
\end{defin}
\begin{rem}\label{convergence}
In order to understand better this notion of convergence, let us
consider a convergent sequence of real numbers $\{x_n\}$. We want to
prove that it is convergent along any ultrafilter $\U$. Let us
consider separately two cases: the first one is when $\U$ is
principal (i.e. there exists $B\subseteq\mathbb N$ such that $\U$ is
the collection of the supersets of $B$. In this case, one says that
$B$ is a basis for $\U$); the second one is when $\U$ is not
principal. In this last case $\U$ is also called \emph{free}. We
need the following classical
\begin{lem}
Let $\U$ be an ultrafilter on a set $X$.
\begin{enumerate}
\item If $\U$ is principal, its basis is a singleton.
\item If $\U$ is not principal, it cannot contain finite sets.
\end{enumerate}
\begin{proof}
\begin{enumerate}
\item Let $B$ the basis for $\U$. If $B$ is not a singleton, we can take a non trivial partition of $B$.
One and only one of the sets of this partition must belong into
$\U$, contradicting the minimality of $B$.
\item Assuming the contrary, let $A\in\U$ be a finite set. Take
$a\in A$. Then one set between $\{a\}$ and $A\setminus\{a\}$ must
belong into $\U$. In the first case $\U$ should be principal with
basis $\{a\}$; in the second one we can repeat the argument until to
obtain a singleton.
\end{enumerate}
\end{proof}
\end{lem}
Coming back to our example, let $\U$ be principal on $\mathbb N$ and
let $\{n_0\}$ be its basis. By definition
$A_{\varepsilon}=\{n\in\mathbb N:|x_n-x_{n_0}|<\varepsilon\}$
contains $n_0$ for all $\varepsilon>0$. Thus $A_{\varepsilon}\in\U$
(by the third property) and consequently $lim_{\U}x_n=x_{n_0}$. On
the other hand, if $\U$ is free, let $x$ be the classical limit of
$\{x_n\}$ and $\varepsilon>0$. One and only one between
$A_{\varepsilon}=\{n:|x_n-x|<\varepsilon\}$ and $\mathbb N\setminus
A_{\varepsilon}$ belongs into $\U$ (by the forth property). But
$\mathbb N\setminus A_{\varepsilon}$ is finite and it follows (by
the lemma) that $A_{\varepsilon}\in\U$ for every $\varepsilon>0$.
\end{rem}
\begin{nota}\label{convergence2}
Notice that we have used that $\{x_n\}$ is a sequence just to
exclude the case $\mathbb N\setminus A_{\varepsilon}\in\U$. A more
refined version of the previous argument however shows that every
bounded net is convergent along a given ultrafilter $\mathcal U$. In
order to prove it one can follow a Bolzano-Weierstrass argument: let
$\{x_a\}_{a\in A}\subseteq[-M,M]$, set $R_1=[-M,0], R_2=(0,M]$ and
$F_i=\{a\in\A:x_a\in R_i\}$. One and only one between $F_1$ and
$F_2$ belongs into the ultrafilter $\mathcal U$ (if it is $F_2$, we
exchange $R_2$ with $\overline{R_2}$ (we find a subset of $A$ which
contains $F_2$ and so it still belongs into $\mathcal U$)). By
repeating this argument, we find a sequence of closed sets $R_n$,
whose diameter halves at each step and containing infinitely many
elements of the net. Now $\bigcap R_n$ is a singleton $\{x\}$ ant it
easy to prove that $lim_{\mathcal U}x_a=x$
\end{nota}
Now we can introduce the notion of ultraproduct. It depends on the
algebraic structure of the objects whose we want to make the
product. Thus there are many kinds of ultraproduct. We are
interested in just two of them: ultraproduct of metric groups and of
type $II_1$ factors. In order to define the ultraproduct of a family
of metric groups we firstly recall what \emph{metric group} means.
\begin{defin}
Let $G$ be a group. A bi-invariant metric on $G$ is a metric on $G$
such that
$$
d(gx,gy)=d(x,y)=d(xg,yg)\,\,\,\,\,\,\,\,\,\,\,\forall x,y,g\in G
$$
\end{defin}
The pair $(G,d)$ is called \emph{metric group}.\\ Similarly one can
define left-invariant or right-invariant metrics, but one can find
examples (see \cite{Pe} Ex.2.1) that show as these concepts are not
good to define the ultraproduct.
\begin{nt}
Let $\{(G_a,d_a)\}_{a\in A}$ be a family of groups equipped with
bi-invariant metrics and $\U$ an ultrafilter on the index set $A$.
We set
$$
G=\{x\in\prod_{a\in A}G_a:sup_{a\in A}d_a(x_a,1_{G_a})<\infty\}
$$
In this way, we assure $lim_{\U}d_a(x_a,1_{G_a})$ exists for any
$x\in G$. So let
$$
N=\{x\in G:lim_{\U}d_a(x_a,1_{G_a})=0\}
$$
We have the following
\end{nt}
\begin{lem}
$N$ is a normal subgroup of $G$.
\begin{proof}
Of course $1_G=\{1_{G_a}\}_{a\in A}\in N$. Let $x,y\in N$, by using
the left invariance and the triangle inequality, one has
$$
d_a(x_ay_a,1_{G_a})=d_a(y_a,x_a^{-1})\leq
d_a(y_a,1_{G_a})+d_a(x_a^{-1},1_{G_a})=d_a(y_a,1_{G_a})+d_a(x_a,1_{G_a})\rightarrow0
$$
Similarly one can prove that if $x\in N$, then also $x^{-1}\in N$.
In order to prove the normality of $N$ we need the hypothesis of
bi-invariance on $d$ (see \cite{Pe} Ex. 2.1). Let $x\in G$ and $n\in
N$. One has
$$
d_a(x_an_ax_a^{-1},1_{G_a})=d_a(x_an_a,x_a)=d_a(n_a,1_{G_a})\rightarrow0
$$
Thus $xnx^{-1}\in N$.
\end{proof}
\end{lem}
Thus the quotient $G/N$ is well-defined as a group and it is easy to
verify it is a metric group with respect to the bi-invariant metric
$$
d(xN,yN)=lim_{\U}d_a(x_a,y_a)
$$
Notice that the metric is well-defined, since $d_a(x_a,y_a)\leq
d_a(x_a,1_{G_a})+d_a(y_a,1_{G_a})$ and thus the net $d_a(x_a,y_a)$
is bounded. Consequently it converges along every ultrafilter (see
Rem.\ref{convergence2}).
\begin{defin}
The metric group $G/N$ is called ultraproduct of the $G_a$'s and it
is denoted by $\prod_{\U} G_a$.
\end{defin}
We will come back to the ultraproduct of metric groups in the next
chapter, when we will describe Radulescu's algebraic approach to the
Conjecture. Now we want to present the construction of the
ultraproduct of type $II_1$ factors $M_a$, which is
\begin{defin}
Let $\{(M_a,tr_a)\}_{a\in A}$ a family of type $II_1$ factors
equipped with normalized traces $tr_a$ and $\U$ an ultrafilter on
$A$. Set
$$
M=\{x\in\prod_{a\in A}M_a:sup||x_a||<\infty\}
$$
and
$$
J=\{x\in M:lim_{\U}tr_a(x_a^*x_a)^{1/2}=0\}
$$
The quotient $M/J$ turns out to be a factor of type $II_1$ with the
trace $tr(x+J)=lim_{\U}tr_a(x_a)$ (but it is not easy to prove! see
\cite{Pe} pg. 18,19 for a sketch, or the original papers by McDuff
(\cite{McD}) and Janssen (\cite{Ja})). It is called
\emph{ultraproduct of the $M_a$'s}. The word \emph{ultrapower} is
referred to the case $M_a=N$ for every $a\in A$.
\end{defin}
The last preliminary notion is a recall of the type $II_1$
hyperfinite factor.
\begin{defin}
A von Neumann algebra $M$ is called approximately finite dimensional
(AFD) if it contains an increasing chain of finite dimensional
subalgebras whose union is strongly dense in $M$.
\end{defin}
It has been already found out by Murray and von Neumann
(\cite{Mu-vN}) that there is substantially a unique (up to von
Neumann algebra isomorphism) AFD factor of type $II_1$, denoted by
$R$. It is called \emph{hyperfinite factor} and it is natural to
expect that it is the smallest type $II_1$ factor, in the sense that
every type $II_1$ factor contains a copy of $R$. Actually, one has
it is the smallest factor of infinite dimension (as Banach space).
One can also describe it explicitly. Let us recall the following
\begin{defin}
Let $G$ be a group and $l^2(G)$ the Hilbert space of all
square-summable complex-valued functions on $G$. Each $g\in G$
defines an operator $\lambda_g:l^2(G)\rightarrow l^2(G)$ in the
following way:
$$
\lambda_g(f)(x)=f(g^{-1}x)
$$
The group von Neumann algebra of $G$, denoted by $VN(G)$, is the
strong operator closure of the subalgebra of $B(l^2(G))$ generated
by all the $\lambda_g$'s.
\end{defin}
\begin{nota}\label{trace}
A group von Neumann algebra is always finite. A trace is determined
by the conditions: $tr(1)=1$ and $tr(g)=0,\forall g\neq1$.
\end{nota}
\begin{rem}\label{immersione}
Notice that $\lambda_g^*=\lambda_{g^{-1}}$. Thus $\lambda_g\in
U(VN(G))$ and the mapping $G\rightarrow U(VN(G))$ defined by
$g\rightarrow\lambda_g$ embeds $G$ into the unitary group of its
group von Neumann algebra.
\end{rem}
\begin{nota}
We recall two classical results on the group von Neumann algebra: it
has separable predual if and only if the group is discrete; it is a
factor if and only if the group is i.c.c., i.e. every conjugacy
class except $\{1_G\}$ is infinite.
\end{nota}
\begin{nota}\label{tensore}
A classical result is that the group von Neumann algebra of
$S_{\infty}^{fin}$ (the group of all the permutations of $\mathbb N$
which fix all but finitely many elements) is the hyperfinite type
$II_1$ factor.\\
Another way to describe the hyperfinite type $II_1$ factor is the
following: $R\cong\bigotimes_{n=1}^{\infty}M_2(\mathbb C)$. Indeed
this von Neumann algebras is a finite factor which contains an
increasing family of factors whose union is strongly dense (by using
$M_2(\mathbb C)\otimes M_2(\mathbb C)\cong M_4(\mathbb C)$). From
this description it follows that $R\otimes R\cong R$, which we will
use later.
\end{nota}
\begin{rem}
The notion of separability for von Neumann algebras cannot be given
with respect to the norm topology, since it is trivial. Indeed, if
$M$ is an infinite dimensional von Neumann algebras, then it
contains a countable family of mutually orthogonal projections, with
which (by using the borelian functional calculus) it is easy to
construct a copy of $l^{\infty}$ into $M$. So the unique von Neumann
algebras which are norm-separable are the finite-dimensional ones.
\end{rem}
The right notion of separability for von Neumann algebras is given
by the following classical
\begin{prop}
Let $M$ be a von Neumann algebras. The following are equivalents:
\begin{enumerate}
\item The predual $M_*$ is norm-separable.
\item $M$ is weakly separable
\item $M$ is faithfully representable into $B(H)$, with $H$
separable.
\end{enumerate}
\end{prop}
\begin{defin}
A von Neumann algebras is called separable if it satisfies one of
the previous conditions.
\end{defin}
After these preliminary notions we are able to enunciate Connes'
embedding conjecture in its original formulation. In order to
simplify notations let us denote $\omega$ a generic free ultrafilter
on $\mathbb N$ and $R^{\omega}$ the ultrapower of $R$ with respect
to $\omega$.
\begin{con}{\bf (A. Connes, \cite{Co})}
Every separable type $II_1$ factor is embeddable into $R^{\omega}$.
\end{con}
\begin{rem}
Assuming Continuum Hypothesis, Ge and Hadwin have proved in
\cite{Ge-Ha} that all the ultrapowers of a fixed $II_1$ factor with
separable predual with regard to a free ultrafilter on the natural
numbers are isomorphic among themselves. More recently, Farah, Hart
and Sherman have proved also the converse: for any separable type
$II_1$ factor $M$ Continuum Hypothesis is \emph{equivalent} to the
statement that all the tracial ultrapowers of $M$ (with regard to a
free ultrafilter on the natural numbers) are isomorphic among
themselves (see \cite{Fa-Ha-Sh}, Th.3.1). On the other hand,
ultrapowers with respect a principal ultrafilter are trivial (being
isomorphic to the factor itself!). It follows that Continuum
Hypothesis together with Connes' embedding conjecture implies the
existence of a universal type $II_1$-factor; \emph{universal} in the
sense that it should contain every type $II_1$ factor. Ozawa have
proved in \cite{Oz2} that such a universal type $II_1$ factor cannot
have separable predual.
\end{rem}
Fortunately we don't have this problem
\begin{prop}\label{separable}
If $\omega$ is non-principal, then $R^{\omega}$ is not separable.
\begin{proof}
We have to prove that $R^{\omega}$ is not faithfully representable
into $B(H)$, with $H$ separable. We recall that if $H$ is separable,
then all the (classical) topologies on $B(H)$ are separable, except
the norm topology (see \cite{Jo}). Moreover, we recall that the
strong topology coincide with Hilbert-Schmidt topology on the
bounded sets. So it is enough to prove that $R^{\omega}$ contains a
non-countable family of unitaries $\{u^{(t)}\}$ such that
$||u^{(t)}-u^{(s)}||_2=\sqrt{2}$ for all $t\neq s$.\\
Let $\{u_n\}\subseteq U(R)$ a sequence of distinct unitaries such
that $u_n\neq1$, for all $n\in\mathbb N$ and $\tau(u_n^*u_m)=0$ for
all $n\neq m$.\\ Let $t\in[\frac{1}{10},1)$, for instance
$t=0,132471...$. Define
$$
I_t=\{1,13,132,1324,13247,132471,...\}
$$
i.e. $I_t$ is the sequence of the approximations of $t$. Clearly,
$\{I_t\}_{t\in[\frac{1}{10},1)}$ is uncountable and $I_t\cap I_s$ is
finite for
all $t\neq s$ (this property forces the choice of $t\geq\frac{1}{10}$!).\\
Now define
$$
u_1^{(t)}=1, u_2^{(t)}=u_2,... u_{12}^{(t)}=u_{12},u_{13}^{(t)}=u_1,
u_{14}^{(t)}=u_2,...u_{131}^{(t)}=u_{131-12},u_{132}^{(t)}=u_1...
$$
i.e. every time we find an element of $I_t$, we start again from
$u_1$. Now define $u^{(t)}=\prod_{n\in\mathbb N}u_n^{(t)}$. Since
$I_t\cap I_s$ is finite (for $t\neq s$), then $u^{(t)}$ and
$u^{(s)}$ have only a finite number of common components. Thus we
have
$$
||u^{(t)}-u^{(s)}||^2_2=lim_{\U}\tau_n((u_n^{(t)}-u_n^{(s)})^*(u_n^{(t)}-u_n^{(s)}))
$$
where $\tau_n$ is the normalized trace on the $n$-th copy of $R$.
Now we observe that
$$
\tau_n((u_n^{(t)}-u_n^{(s)})^*(u_n^{(t)}-u_n^{(s)}))=\left\{
                                       \begin{array}{ll}
                                         0\,\,\,\,\,\,\,\,\,\,\,if\,\,\,\,\,u_n^{(t)}=u_n^{(s)} \\
                                         2\,\,\,\,\,\,\,\,\,\,\,if\,\,\,\,\,u_n^t\neq
u_n^s
                                       \end{array}
                                     \right.
$$
Since $u_n^{(t)}=u_n^{(s)}$ only on a finite set and since $\U$ is
free (and thus it does not contain finite sets), it follows that
$$
lim_{\U}\tau_n((u_n^{(t)}-u_n^{(s)})^*(u_n^{(t)}-u_n^{(s)}))=2
$$
and thus $||u^{(t)}-u^{(s)}||_2=\sqrt{2}$.
\end{proof}
\end{prop}
\begin{nota}
Non-separability of $R^{\omega}$ has been already proved by several
authors (\cite{Fe} and, in greater generality, \cite{Po} Prop. 4.3).
We have preferred this proof because it is constructive in the sense
that will be more clear in the following section: it allows
(together with th. \ref{product}) to produce examples of uncountable
groups which embed trace-preserving into $U(R^{\omega})$ and to
generalize a theorem by R\u adulescu (see also \cite{Ca-Pa}).
\end{nota}
\section{The algebraic approach}
The idea of the algebraic approach is attaching the following weaker
version of Connes' embedding conjecture.
\begin{con}{\bf (Connes' embedding conjecture for groups)}
For every countable i.c.c. group $G$, the group von Neumann algebra
$VN(G)$ embeds into a suitable ultrapower $R^{\omega}$.
\end{con}
R\u adulescu in \cite{Ra1} has worked to find a characterization for
those groups which satisfy this weaker version of the Conjecture.
\begin{defin}
Let $U(n)$ be the unitary group of order $n$, i.e. the group of
$n\times n$ matrices with complex entries and such that
$u^*u=uu^*=1$. The normalized Hilbert-Schmidt distance on $U(n)$ is
$$
d_{HS}(u,v)=||u-v||_2=\sqrt{\frac{1}{n}\sum_{i,j=1}^n|u_{ij}-v_{ij}|}=\frac{1}{\sqrt{n}}\sqrt{tr((u-v)(u-v)^*)}
$$
\end{defin}
Bi-invariance of this metric follows from the main property of the
trace: $tr(ab)=tr(ba)$.
\begin{defin}
A group $G$ is called hyperlinear if it embeds into a suitable
ultraproduct of unitary groups of finite rank (equipped with
Hilbert-Schmidt distance).
\end{defin}
Our purpose is to prove the following characterization theorem
\begin{teo}{\bf (R\u adulescu)}\label{algebrico}
The following conditions are equivalent
\begin{enumerate}
\item Connes' embedding conjecture for groups is true.
\item Every countable i.c.c. group is hyperlinear.
\end{enumerate}
\end{teo}
This theorem was firstly proved by R\u adulescu in the countable
case. L. Paunescu and the present author have generalized it to the
continuous one. We present our result later (see
Cor.\ref{continuous}). Now we need some preliminary results.
\begin{lem}\label{uno}
Let $M=\prod_{\U} M_a$ be a type $II_1$ factor obtained as
ultraproduct of type $II_1$ factors $M_a$, equipped with normalized
traces $tr_a$, with regard to an ultrafilter $\U$ on the index set.
Then
$$
U(M)=\prod_{\U}U(M_a)
$$
i.e. the unitary group of the ultraproduct is the ultraproduct of
the unitary groups.
\begin{proof}
The inclusion $\supseteq$ is obvious, since the multiplication in
the ultraproduct is pointwise. Conversely, let $v_a\in M_a$ such
that $v=\prod_{\U}v_a$ is unitary, i.e. $\prod_{\U}v_a^*v_a=1$. We
have to prove that there exist unitaries $u_a$ such that
$\prod_{\U}u_a=\prod_{\U}v_a$. Let $v_a=u_a|v_a|$ the polar
decomposition of $v_a$. Since $M_a$ is a type $II_1$ factor, we can
extend the partial isometry $u_a$ to a unitary operator. So we can
assume that $u_a$ is unitary. Now we can verify that they are just
the unitaries which we are looking for. Indeed
$$
\prod_{\U}v_a=\prod_{\U}u_a|v_a|=\prod_{\U}u_a\prod_{\U}|v_a|=\prod_{\U}u_a\prod_{\U}(v_a^*v_a)^{1/2}=\prod_{\U}u_a(\prod_{\U}v_a^*v_a)^{1/2}=\prod_{\U}u_a
$$
\end{proof}
\end{lem}
\begin{prop}{\bf (Elek-Szab\'{o}, \cite{El-Sz})}\label{perdue}
Let $G$ be a group such that for any finite $F\subseteq G$ and any
$\varepsilon>0$ there exist a natural number $n$ and a map
$\theta:F\rightarrow U(n)$ such that
\begin{enumerate}
\item if $g,h,gh\in F$, then
$||\theta(gh)-\theta(g)\theta(h)||_2<\varepsilon$
\item if $1_G\in F$, then $||\theta(1_g)-1_{U(n)}||_2<\varepsilon$
\item for all distinct $x,y\in F$, $||\theta(x)-\theta(y)||\geq1/4$
\end{enumerate}
Then $G$ is hyperlinear.
\begin{proof}
Choosing $\varepsilon=1/n$, we have a family of maps
$\theta_{F,1/n}:F\rightarrow U(F,n)$. We set
$$
A=\{(F,1/n), F\subseteq G\,\,\, finite, n\geq1\}
$$
partially ordered in a natural way. Let $\U$ be a free ultrafilter
on $A$ containing every subset of the form $\{(H,1/m):H\supseteq
F,m\geq n\}$. Now we consider the map
$$
\theta: G\ni
g\rightarrow\prod_{\U}\theta_{F,1/n}(g)\in\prod_{\U}U(F,n)
$$
We have to prove that this map is a monomorphism. Let $d_{HS}^{F,n}$
and $d_{HS}$ respectively the Hilbert-Schmidt distance on $U(F,n)$
and on the ultraproduct $\prod_{\U}U(F,n)$. We have
$$
d_{HS}(\theta(hg)-\theta(h)\theta(g))=lim_{\U}d_{HS}^{F,n}(\theta_{F,1/n}(gh),\theta_{F,1/n}(g)\theta_{F,1/n}(h))\leq
$$
$$
\leq lim_{\U}1/n
$$
Now we use the particular choice of the ultrafilter in order to
conclude that the previous limit must be zero. In a similar way (by
using the second property) one can easily prove that $\theta$ is
unital. Thus it is an homomorphism. Injectivity follows from the
third property applied to a similar argument.
\end{proof}
\end{prop}
\begin{nota}
Also the converse of the previous proposition is true (see
\cite{El-Sz}). Moreover, this proposition shows that the notion of
hyperlinearity does not depend on the choice of the ultrafilter.
\end{nota}
\begin{rem}
The previous proposition can be viewed in the following way: if one
can approximate every finite subset of $G$ with a unitary group of
finite rank, then $G$ is hyperlinear. A fundamental application of
this fact is the following
\end{rem}
\begin{cor}\label{due}
Let $(G,d)$ be a metric group containing an increasing chain of
subgroups isomorphic to $U(n), n\in\mathbb N$, whose union is dense
in $G$ and such that $d|_{U(n)}=d_{HS}$. Then a group $H$ is
hyperlinear if and only if it embeds into a suitable ultrapower of
$G$.
\begin{proof}
If $H$ is hyperlinear, then it embeds into a suitable ultraproduct
of unitary groups. This ultraproduct, by hypothesis, embeds into the
same ultraproduct of $G$. Conversely, let $\theta$ be the embedding
of $H$ into $\prod_{\U}G$. Let $F=\{f_1,...f_k\}\subseteq H$ and
$\varepsilon>0$. Let $m$ be a natural number and $(\theta(f_i))_m$
the m-th component of $\theta(f_i)\in G^{\omega}$. Since $F$ is
finite and $\theta$ is an embedding, we can choose $m$ such that
$(\theta(f_i))_m$ are all distinct (take $m$ in the intersection
among the sets on which the $\theta(f_i)_m$'s differ). This
intersection cannot be empty, since it is finite intersection of
sets belonging into $\omega$). By hypothesis, there exist $u_i\in
U(n_i)$ such that $d_{HS}((\theta(f_i))_m-u_i)<\varepsilon$. Let
$n=max\{n_i,i=1,...k\}$, we can identify $u_i\in U(n)$. Moreover we
can assume that $u_i$ are all distinct. So, we can define
$\theta_{F,\varepsilon}(f_i)=u_i$: the first two properties are
clearly satisfied; the third one is not obvious and one has to
follow a trick known as \emph{amplification} (see \cite{El-Sz} and
\cite{Ra1}).
\end{proof}
\end{cor}
\begin{lem}\label{tre}
$$
R^{\omega}\otimes R^{\omega}\subseteq R^{\omega}
$$
\begin{proof}
At first we observe that $(R\otimes R)^{\omega}\cong R^{\omega}$, by
using the isomorphism $(x_n\otimes
y_n)_n\rightarrow(\theta(x_n\otimes y_n))_n$, where $\theta$ is an
isomorphism between $R\otimes R$ and $R$. It remains to embed
$R^{\omega}\otimes R^{\omega}$ into $(R\otimes R)^{\omega}$. We can
do this by using the embedding $(x_n)_n\otimes(y_n)_n\rightarrow
(x_n\otimes y_n)_n$.
\end{proof}
\end{lem}
\begin{lem}\label{quattro}
Let $G$ be an i.c.c. group and $\theta:G\rightarrow U(M)$ a unitary
faithfully representation on a finite von Neumann algebra $M$ with
trace $\tau$. Then $|\tau(\theta(g))|<1, \forall g\in G,g\neq1_G$.
\begin{proof}
Certainly $|\tau(\theta(g))|\leq1$, since the unitary elements have
trace in absolute value $\leq1$. So we assume
$\tau(\theta(g_0))=\lambda$, with $\lambda$ a complex number with
norm one, and we prove that $g_0=1_G$. Take
$u=\lambda^*\theta(g_0)$. This unitary element has trace one and
thus it must be the identity. Indeed
$$
\tau((u-1)^*(u-1))=2-2Re(\tau(u))=0
$$
So $\theta(g_0)=\lambda1$ and thus, setting $h=gg_0g^{-1}$, we have
$\theta(h)=\lambda$ and thus $h=g_0$ (by the faithfulness).
Therefore $g_0$ is the unique element of its conjugacy class and
then $g_0=1_G$.
\end{proof}
\end{lem}
Now we can prove Radulescu's characterization theorem.
\begin{proof}
We have to prove that $VN(G)$ embeds into $R^{\omega}$ if and only
if $G$ is hyperlinear. We start assuming that $VN(G)$ embeds into
$R^{\omega}$. Recalling rem.\ref{immersione} we have that $G$ embeds
into $U(VN(G))$ and then into $U(R^{\omega})$. It follows that $G$
embeds into $\prod_{\omega} U(R)$ (by the lemma \ref{uno}). Now we
recall that $R$ contains an increasing family of weakly dense finite
dimensional von Neumann factors. Thus $U(R)$ contains an increasing
family of subgroups isomorphic to $U(n)$ whose union is dense in
$U(R)$ (the density follows from the normality of the trace). Of
course we have the restriction property of the distance. So we can
use Cor.\ref{due} to conclude that $G$ must be hyperlinear.
Conversely, by the hypothesis and by Cor. \ref{due}, $G$ embeds into
$U(R^{\omega})$. Let $\theta_1$ be such an embedding and $g\in G$,
$g\neq1$. By Lemma \ref{quattro}, we have $|\tau(\theta_1(g)|<1$ We
define a new embedding $\theta_2=\theta_1\otimes\theta_1$. This is
still an embedding into $U(R^{\omega})$, by Lemma \ref{tre}.
Moreover $\tau(\theta_2(g))=\tau(\theta_1(g))^2$. By induction we
can construct a sequence of embedding
$\theta_n=\theta_{n-1}\otimes\theta_1$ and we have
$lim|\tau(\theta_n(g))|=0$ for each $g\neq1$. Now, since $G$ is
countable, we can write $G=\bigcup_{i=1}^{\infty} A_i$, with
$A_1\subseteq A_2\subseteq ...$ are all finite subsets of $G$. For
each $k\in\mathbb N$ choose $\theta_{n_k}$ such that
$|\tau(\theta_{n_k}(g))|<2^{-k}$ for all $g\in A_k$, $g\neq1$. Let
us denote $\lambda_k=\theta_{n_k}$. Moreover, if $x\in
U(R^{\omega})$, $x_n$ denotes the $n$-th component of $x$. Lastly
$\tau_n$ denotes the trace on the $n$-th copy of $R$. With these
notations, we have
$$
lim_{\omega}|\tau_n(\lambda_k(g)_n)|=|\tau(\lambda_k(g))|<2^{-k}
$$
Thus, there exists $m_k$ such that
$\tau_{m_k}(\lambda_k(g)_{m_k})<2^{-k}$. We define
$\pi_{m_k}:G\rightarrow U(R)$ by setting
$\pi_{m_k}(g)=\lambda_k(g)_{m_k}$ and lastly
$\pi=\prod_{\omega}\pi_{m_k}:G\rightarrow U(R^{\omega})$ by setting
$\pi(g)=\prod_{\omega}\pi_{m_k}(g)$. It is still an embedding and
verifies the fundamental property that $\tau(\pi(g))=0$ for all
$g\neq1$ and $\tau(\pi(1))=1$, indeed for $g\neq1$
$$
|\tau(\pi(g))|=lim_{\omega}|\tau_k(\pi(g)_k)|=lim_{\omega}|\tau_k(\lambda_k(g)_{m_k})|=lim_{\omega}|\tau_{m_k}(\lambda_k(g)_{m_k})|<lim_{\omega}2^{-k}=0
$$
So the trace of $\pi(g)$ is exactly the trace of $g$, viewed into
the group von Neumann algebra (see Note \ref{trace}). This means
just that we can extend this embedding to $VN(G)$ and find an
identification between $VN(G)$ and a subalgebra of $R^{\omega}$.
\end{proof}
The previous proof is quite technical and strongly depending on the
hypothesis of countability of $G$. We can simplify and extend it to
the uncountable case by using a concept of product between
ultrafilters. This notion, together with Prop. \ref{separable}, also
allows to prove that the von Neumann algebra of the free group on a
continuous family of generators $VN(\mathbb F_{\aleph_c})$ is
embeddable into $R^{\omega}$. More applications of the product
between ultrafilters can
be found in \cite{Ca-Pa}.\\\\
In order to prove it, we recall the classical notion of tensor
product between ultrafilters and we prove that
$(R^{\omega})^{\omega'}\cong R^{\omega\otimes\omega'}$.
\begin{defin}
Let $\omega,\omega'$ be two ultrafilters on $\mathbb N$. We define
$$
B\in\omega\otimes\omega'\Leftrightarrow\{k\in\mathbb N:\{n\in\mathbb
N:(k,n)\in B\}\in\omega'\}\in\omega
$$
\end{defin}
\begin{teo}\label{product}
$$
(R^{\omega})^{\omega'}\cong R^{\omega\otimes\omega'}
$$
\begin{proof}
Since operations are component-wise we don't have algebraic problem.
We have only to prove that those factors have the same trace. So we
have to prove that
$$
lim_{k\rightarrow\omega}lim_{n\rightarrow\omega'}x_n^k=lim_{(k,n)\rightarrow\omega\times\omega'}x_n^k
$$
Let $x=lim_{k\rightarrow\omega}lim_{n\rightarrow\omega'}x_n^k$.
Fixed $\varepsilon>0$, we set
$$
A=\{k\in\mathbb
N:|lim_{n\rightarrow\omega'}x_n^k-x|<\frac{\varepsilon}{2}\}\in\omega
$$
and
$$
A_k=\{n\in\mathbb
N|x_n^k-lim_{n\rightarrow\omega'}x_n^k|<\frac{\varepsilon}{2}\}\in\omega'
$$
So
$$
B=\{(k,n)\in\mathbb N^2:k\in A,n\in A_k\}\subseteq\{(k,n)\in\mathbb
N^2:|x_n^k-x|<\varepsilon\}
$$
Since $B\in\omega\otimes\omega'$ the proof is complete.
\end{proof}
\end{teo}
In the proof of Radulescu's theorem we have used the fact that a
countable group is hyperlinear if and only if it embeds into
$U(R^{\omega})$. The \emph{only if} part is no longer true for
uncountable groups because they can be too big. So, the right way to
extend Radulescu's theorem to the uncountable case is the following
\begin{cor}\label{continuous}
For an i.c.c. group $G$, the following statements are equivalent
\begin{enumerate}
\item $G$ (not necessarily countable) embeds into $U(R^{\omega})$.
\item $VN(G)$ is embeddable into $R^{\omega}$.
\end{enumerate}
\begin{proof}
Radulescu's proof of the implication $2.\Rightarrow1.$ does not
depend on the countability of $G$. Conversely, we can follow
Radulescu's proof and define $\theta_n$. Then, we define
$\theta(g)=\{\theta_n(g)\}_{n\in\mathbb N}$. It is an embedding into
$U((R^{\omega})^{\omega})$. By Th.\ref{product} one can look at
$\theta$ as an embedding into $U(R^{\omega\times\omega})$. Now
$\tau(\theta(g))=lim_{\omega}\tau(\theta_n(g))=0$, whenever
$g\neq1$.
\end{proof}
\end{cor}
Here is a nice application of the product between ultrafilters and
of the construction of Prop.\ref{separable}. Let $F_{\aleph_c}$ be
the free group on a continuous family of generators.
\begin{cor}\label{uncountable}
$VN(F_{\aleph_c})$is embeddable into $R^{\omega}$.
\begin{proof}
It is enough to prove that $F_{\aleph_c}$ is embeddable into
$U(R^{\omega})$ and that such an embedding $\theta$ preserves the
trace, i.e. $\tau(\theta(g))=0$ if $g\neq1$ and $\tau(\theta(1))=1$.
Since $F_{\infty}$ (free group countably generated) is hyperlinear,
we have a sequence $\{u_n\}\subseteq U(R^{\omega})$ such that
\begin{enumerate}
\item $\tau(u_n)=0$ for all $n\in\mathbb N$
\item $\tau(u_n^*u_m)=0$ for all $n\neq m$
\item $u_n$ have no relations between themselves
\end{enumerate}
This sequence is simply the image of the generators of $F_{\infty}$
into $U(R^{\omega})$. Now, we apply the construction of the proof of
Prop.\ref{separable}. By using Th.\ref{product}, we find a copy of
$F_{\aleph_c}$ into $U(R^{\omega})$ such that the desired property
on the trace is satisfied.
\end{proof}
\end{cor}
\begin{nota}
In this last note we want to describe briefly the actual situation
of the research around hyperlinear groups. Indeed, in the hope that
Connes' embedding conjecture is true, one can try to prove that any
group is hyperlinear. This problem is still open, but there are some
positive partial results.
\begin{enumerate}
\item We recall the following
\begin{defin}
A group is called residually finite if the intersection of all its
normal subgroups of finite index is trivial.
\end{defin}
Clearly, finite groups are residually finite. A non-trivial example
of residually finite groups is given by the free groups (\cite{Sa}).
Anyway\\\\
Every residually finite group is hyperlinear (see \cite{Pe}, ex.
4.2)
\item We recall the following
\begin{defin}
A group $G$ is called amenable if for any finite $F\subseteq G$ and
$\varepsilon>0$, there exists a finite $\Phi\subseteq G$ such that
for every $g\in F$,
$$
|g\Phi\Delta\Phi|<\varepsilon|\Phi|\,\,\,\,\,\,\,\,\,\,\,\,\,\,\,\,\,\,
Folner\,\,\,condition
$$
where $\Delta$ stands for the symmetric difference: $A\Delta
B=(A\cup B)\setminus (A\cap B)$.
\end{defin}
For instance, compact groups are amenable (One should use an
equivalent definition of amenability, linked to the measure theory.
Then
amenability of compact groups follows from the finiteness of the Haar measure).\\\\
Every amenable group is hyperlinear (see \cite{Pe}, ex. 4.4).\\\\
Gromov introduced in \cite{Gr} the notion of initially subamenable
groups: these are groups for which every finite subset can be
multiplicatively embedded into an amenable group. By
Prop.\ref{perdue} (and its converse) it follows that hyperlinearity
is a local property. Thus, also initially subamenable groups are
hyperlinear.
\end{enumerate}
The class of initially subamenable groups is the largest among those
we have a result about hyperlinearity (residually finite groups are
initially subamenable (see the next proposition)). Thom has been the
first who had found an example of hyperlinear group which is not
initially subamenable (see \cite{Th}).
\end{nota}
Let us conclude this section with the proof that every residually
finite group is initially subamenable. This fact is well known but
it seems to impossible to find references.
\begin{defin}
A group is called initially subfinite if every finite subset can be
multiplicatively embedded into a finite group.
\end{defin}
Since finite groups are amenable, it is enough to prove the
following
\begin{prop}\label{liviu}
Every residually finite group is initially subfinite.
\begin{proof}
Let $G$ be a residually finite group and $F\subseteq G$ be a finite
subset. We set $\overline F=F\cup\{xy, x,y\in F\cup F^{-1}\}$. So
$\overline F$ is still finite. Now, for each $x\in\overline F$ there
exists a normal subgroup $G_x\unlhd G$ with finite index not
containing $x$ (up to the case $x=1$). Let $H=\bigcap_{x\in\overline
F} G_x$. So $\overline F\cap H$ is the empty set or the identity.
Moreover $H$ is still a normal subgroup of $G$ with finite index
(since it is finite intersection of normal subgroups with finite
index). Let $\pi:G\rightarrow G/H$ the canonical projection. It
remains to observe that $\pi|_{\overline F}$ is an embedding which
preserves the multiplication from $\overline F$ into the finite
group $G/H$.
\end{proof}
\end{prop}
The converse of the previous proposition is false: in \cite{Th} one
can find an example of initially subfinite group which is not
residually finite.
\section{The topological approach}
In this section we want to describe the topological approach by
Haagerup and Winslow (see \cite{Ha-Wi1} and \cite{Ha-Wi2}). Let $H$
be a Hilbert space and $vN(H)$ the set of von Neumann algebras
acting on $H$, this topological approach is based on the definition
of a topology on $vN(H)$, named Effros-Marechal topology. Indeed
they were the firsts who have introduced this topology and have
studied its properties (see \cite{Ef} and \cite{Ma}); but merely
Haagerup and Winslow, thirty years later, have argued the link
between this topology and Connes' embedding conjecture.
\\
There are three different ways to describe the Effros-Marechal
topology and one can find  in \cite{Ha-Wi1} (th.2.8) the proof that
these ways are truly equivalent. Here we shall give only the
definitions and we shall describe some interesting properties
without giving proofs. Here is the first definition
\begin{defin}
The Effros-Marechal topology on $vN(H)$ is the weakest topology such
that for every $\phi\in B(H)_*$ the mapping
$$
vN(H)\ni M\rightarrow||\phi|_M||
$$
is continuous.
\end{defin}
The second definition of the Effros-Marechal topology come from a
more general definition by Effros (see \cite{Ef2})
\begin{nt}
Let $X$ be a compact Hausdorff space, $c(X)$ the set of closed
subsets of $X$ and $\omega(x)$ the set of the neighborhoods of a
point $x\in X$. Let $\{C_a\}\subseteq c(X)$ and
$$
\underline{lim}C_a=\{x\in X:\forall U\in\omega(x), U\cap
C_a\neq\emptyset\,\,eventually\}
$$
$$
\overline{lim}C_a=\{x\in X:\forall U\in\omega(x), U\cap
C_a\neq\emptyset\,\,frequently\}
$$
\end{nt}
Effros has proved that there is only a topology on $c(X)$, whose
convergence is described by the conditions
$$
C_a\rightarrow
C\,\,\,\,\,\,\,\,\,\,\,\,\,\,\,iff\,\,\,\,\,\,\,\,\,\,\,\,\,\,\,\overline{lim}C_a=\underline{lim}C_a=C
$$
Since the unit ball $Ball(M)$ of a von Neumann algebra $M$ is weakly
compact, one can use this notion of convergence in our setting.
\begin{defin}\label{second}
Let $\{M_a\}\subseteq vN(H)$ be a net. The Effros-Marechal topology
is described by the following notion of convergence:
$$
M_a\rightarrow M
\,\,\,\,\,\,\,\,\,\,\,\,\,\,\,iff\,\,\,\,\,\,\,\,\,\,\,\,\,\,\,\overline{lim}Ball(M_a)=\underline{lim}Ball(M_a)=Ball(M)
$$
\end{defin}
The third definition is by introducing a further notion of
convergence in $vN(H)$. First of all we need some definitions.
\begin{nt}
Let $x\in B(H)$, $so^*(x)$ denotes the set of the neighborhoods of
$x$ with respect to the strong* topology.
\end{nt}
\begin{defin}
Let $\{M_a\}\subseteq vN(H)$ be a net. We set
$$
liminfM_a=\{x\in B(H):\forall U\in so^*(x), U\cap
M_a\neq\emptyset\,\, eventually\}
$$
\end{defin}
By th. 2.6 in \cite{Ha-Wi1}, $liminfM_a$ can be thought as the
largest element in $vN(H)$ whose unit ball is contained in
$\underline{lim}Ball(M_a)$. This suggests to define $limsupM_a$ as
the smallest element in $vN(H)$ whose unit ball contains
$\overline{lim}Ball(M_a)$, that is clearly
$(\overline{lim}Ball(M_a))''$. So we have quite naturally the
following
\begin{defin}
Let $\{M_a\}\subseteq vN(H)$ be a net. We set
$$
limsupM_a=(\overline{lim}Ball(M_a))''
$$
\end{defin}
Now here is the third description of the Effros-Marechal topology
\begin{defin}
The Effros-Marechal topology on $vN(H)$ is described by the
following notion of convergence:
$$
M_a\rightarrow
M\,\,\,\,\,\,\,\,\,\,\,\,\,\,\,\,iff\,\,\,\,\,\,\,\,\,\,\,\,\,\,\,\,liminf
M_a=limsupM_a=M
$$
\end{defin}
We recall that in \cite{Ha-Wi1}, th. 2.8, they have shown that these
three definition of the Effros-Marechal topology are equivalent.
\\\\
Connes' embedding conjecture regards the behaviour of separable type
$II_1$ factors. So we are interested in the case in which $H$ is
separable. In this case it happens that the Effros-Marechal topology
is metrizable, second countable and complete (i.e. $vN(H)$ is a
Polish space). Moreover, a possible distance is given by the
Hausdorff distance between the unit balls:
$$
d(M,N)=max\{sup_{x\in Ball(M)}\{inf_{y\in Ball(N)}d(x,y)\},sup_{x\in
Ball(N)}\{inf_{y\in Ball(M)}d(x,y)\}\}
$$
where $d$ is a metric on the unit ball of $B(H)$ which induces the
weak topology (remember that the weak topology on
$Ball(B(H))$ is metrizable whenever $H$ is separable).\\
There are many interesting results about the Effros-Marechal
topology in the case of separability of $H$. For example, the sets
of factors of each of the types $I_n,n\in\mathbb N$, $II_1$,
$II_{\infty}$, $III_{\lambda},\lambda\in[0,1]$ are Borel subsets of
$vN(H)$ without being $G_{\delta}$-sets, or many others (see
\cite{Ha-Wi1} sections 4. and 5.). Anyway it is in the second paper
\cite{Ha-Wi2} that Haagerup and Winslow have begun studying density
problems relative to the Effros-Marechal topology. What important
subsets of $vN(H)$ are dense? They have found a lot of interesting
results and, above all, an equivalent condition to Connes' embedding
conjecture.
\begin{rem}\label{strong}
In the case in which $H$ is separable and $\{M_a\}=\{M_n\}$ is a
sequence, the definition of the Effros-Marechal topology may be
simplified by using the third definition. In particular we have
$$
liminf M_n=\{x\in B(H): \exists\{x_n\}\in\prod
M_n\,\,s.\,\,t.\,\,x_n\rightarrow^{s^*}x\}
$$
\end{rem}
\begin{nt}
$$
\Im_{I_{fin}}\,\,is\,\,the\,\,set\,\,of\,\,type\,\,I\,\,finite\,\,factors\,\,acting\,\,on\,\,H
$$
$$
\Im_{I}\,\,is\,\,the\,\,set\,\,of\,\,type\,\,I\,\,factors\,\,acting\,\,on\,\,H
$$
$$
\Im_{AFD}\,\,is\,\,the\,\,set\,\,of\,\,approximately\,\,finite\,\,dimensional\,\,factors\,\,acting\,\,on\,\,H
$$
$$
\Im_{inj}\,\,is\,\,the\,\,set\,\,of\,\,injective\,\,factors\,\,acting\,\,on\,\,H
$$
We recall that a von Neumann algebra $M\subseteq B(H)$ is called
injective if it is the range of a projection of norm $1$. For
example every type I von Neumann algebra is injective.
\end{nt}
\begin{teo}\label{haa}
The following statements are equivalent:
\begin{enumerate}
\item $\Im_{I_{fin}}$ is dense in $vN(H)$
\item $\Im_I$ is dense in $vN(H)$
\item $\Im_{AFD}$ is dense in $vN(H)$
\item $\Im_{inj}$ is dense in $vN(H)$
\item Connes' embedding conjecture is true
\end{enumerate}
Moreover, a separable type $II_1$ factor $M$ is embeddable into
$R^{\omega}$ if and only if $M\in\overline{\Im_{inj}}$.
\begin{proof}
As $\Im_{I_{fin}}\subseteq\Im_I\subseteq\Im_{AFD}$, the implications
$1.\Rightarrow2.\Rightarrow 3.$ are trivial. The implication
$3.\Rightarrow1.$ follows from the fact that AFD factors contain an
increasing chain of type $I_{fin}$ factors, whose union is weakly
dense and from the second definition of the Effros-Marechal topology
(Def.\ref{second}). The equivalence between 3. and 4. is a theorem
by A. Connes (\cite{Co}). The equivalence between 4. and 5. is the
theorem by Haagerup and Winslow (\cite{Ha-Wi2}, Cor.5.9).\\
Also the last sentence is proved in \cite{Ha-Wi2} (see Th. 5.8).
\end{proof}
\end{teo}
Now we want to give a sketch of Haagerup-Winslow's proof of a
theorem by Kirchberg, which gives probably the most unexpected
equivalent condition to Connes' embedding conjecture. Let us recall
some concepts on the tensor product of $C^*$-algebras. A complete
introduction can be found in the forth chapter of the first book by
Takesaki (\cite{Ta1}).
\begin{rem}
The algebraic tensor product of two $C^*$-algebras is a *algebra in
a natural way, by setting
$$
(x_1\otimes x_2)(y_1\otimes y_2)=x_1x_2\otimes y_1y_2
$$
$$
(x_1\otimes x_2)^*=x_1^*\otimes x_2^*
$$
Nevertheless it is not clear how one can define a norm to obtain a
$C^*$-algebra (Notice that the product of the norms is not in
general a norm on the algebraic tensor product).
\end{rem}
\begin{defin}
Let $A_1,A_2$ be two $C^*$-algebras and $A_1\otimes A_2$ their
algebraic tensor product. A norm $||\cdot||_{\beta}$ on $A_1\otimes
A_2$ is called $C^*$-norm if the followings hold
\begin{enumerate}
\item $||xy||_{\beta}\leq||x||_{\beta}||y||_{\beta}$, for all
$x,y\in A_1\otimes A_2$
\item $||x^*x||_{\beta}=||x||_{\beta}^2$, for all $x\in A_1\otimes
A_2$
\end{enumerate}
If $||\cdot||_{\beta}$ is a $C^*$-norm on $A_1\otimes A_2$,
$A_1\otimes_{\beta}A_2$ stands for the completion of $A_1\otimes
A_2$ with respect to $||\cdot||_{\beta}$. It is a $C^*$-algebra.
\end{defin}
Unlucky there is no a unique $C^*$-norm on $A_1\otimes A_2$ in
general, but one can construct by hands at least two of them.
\begin{defin}
$$
||x||_{max}=sup\{||\pi(x)||,
\pi\,\,*representation\,\,of\,\,the\,\,*algebra\,\,A_1\otimes A_2\}
$$
This norm is called projective or Turumaru's norm (\cite{Tu}). One
can prove that this norm is a $C^*$-norm and the completion of
$A_1\otimes A_2$ with respect to it is denoted by
$A_1\otimes_{max}A_2$.
\end{defin}
The projective norm has the following universal property (see
\cite{Ta1}, IV.4.7)
\begin{prop}\label{universal}
Given $C^*$-algebras $A_1,A_2,B$. If $\pi_i:A_i\rightarrow B$ are
homomorphisms with commuting ranges, then there exists a unique
homomorphism $\pi:A_1\otimes_{max}A_2\rightarrow B$ such that
\begin{enumerate}
\item $\pi(x_1\otimes x_2)=\pi_1(x_1)\pi_2(x_2)$
\item $\pi(A_1\otimes_{max}A_2)=C^*(\pi_1(A_1),\pi_2(A_2))$
\end{enumerate}
\end{prop}
\begin{defin}
$$
||x||_{min}=sup\{||(\pi_1\otimes\pi_2)(x)||,
\pi_1\,\,representation\,\,of\,\,A_1,
\pi_2\,\,representation\,\,of\,\,A_2\}
$$
This norm is called injective or Guichardet's norm (\cite{Gu}). One
can prove that this norm is a $C^*$-norm and the completion of
$A_1\otimes A_2$ with respect to it is denoted by
$A_1\otimes_{min}A_2$.
\end{defin}
\begin{rem}
Clearly $||\cdot||_{min}\leq||\cdot||_{max}$, since representations
of the form $\pi_1\otimes\pi_2$ are particular *representation of
the algebraic tensor product $A_1\otimes A_2$. These norms are
different, in general, as Takesaki has shown in \cite{Ta2}. More
recently Junge and Pisier have shown, in \cite{Ju-Pi}, that
$B(l^2)\otimes_{min}B(l^2)\neq B(l^2)\otimes_{max}B(l^2)$. Notation
$||\cdot||_{max}$ reflects the obvious fact that there are no
$C^*$-norm greater than that one. Notation $||\cdot||_{min}$ has the
same justification, but it is harder to prove:
\begin{teo}{\bf (Takesaki, \cite{Ta2})}
$||\cdot||_{min}$ is the smallest $C^*$-norm among those on
$A_1\otimes A_2$.
\end{teo}
\end{rem}
\begin{defin}
Let $G$ be a locally compact group. By using the Haar measure, one
can consider $L^1(G)$. The universal $C^*$-algebra of $G$ is the
envelopping $C^*$-algebra of $L^1(G)$, i.e. the completion of
$L^1(G)$ with respect to the norm $||f||=sup_{\pi}||\pi(f)||$, where
$\pi$ runs over all non-degenerate *representation of $L^1(G)$ in a
Hilbert space. This norm makes sense by virtue of the classical
result: a *homomorphism of an involutive Banach algebra into a
$C^*$-algebra is contractive.
\end{defin}
\begin{rem}
Let $\mathbb F_{\infty}$ the free group countably generated. It is a
locally compact group with respect to the discrete topology, so we
can consider its universal $C^*$-algebra, $C^*(\mathbb F_{\infty})$.
$\mathbb F_{\infty}$ can be canonically embedded into $U(C^*(\mathbb
F_{\infty}))$. Unitaries corresponding to such an embedding are
called universal.
\end{rem}
Here is Kirchberg's theorem (\cite{Ki}).
\begin{teo}\label{kirchi}
The following statements are equivalent
\begin{enumerate}
\item $C^*(\mathbb F_{\infty})\otimes_{min}C^*(\mathbb F_{\infty})=C^*(\mathbb F_{\infty})\otimes_{max}C^*(\mathbb F_{\infty})$
\item Connes' embedding conjecture is true.
\end{enumerate}
\begin{proof}
By using Th.\ref{haa} we can prove the following implications:
\begin{enumerate}
\item If $\Im_{I_{fin}}$ is dense in $vN(H)$, then $C^*(\mathbb F_{\infty})\otimes_{min}C^*(\mathbb F_{\infty})=C^*(\mathbb F_{\infty})\otimes_{max}C^*(\mathbb F_{\infty})$
\item If
$C^*(\mathbb F_{\infty})\otimes_{min}C^*(\mathbb
F_{\infty})=C^*(\mathbb F_{\infty})\otimes_{max}C^*(\mathbb
F_{\infty})$, then $\Im_{inj}$ is dense in $vN(H)$.
\end{enumerate}
(Proof of 1.)\\ Let $\pi$ be a
*representation of the algebraic tensor product
$C^*(\mathbb F_{\infty})\otimes C^*(\mathbb F_{\infty})$ into
$B(H)$. Since $C^*(\mathbb F_{\infty})$ is separable, we can assume
that $H$ is separable. In this way
$$
A=\pi(C^*(\mathbb F_{\infty})\otimes\mathbb
C1)\,\,\,\,\,\,\,\,\,\,\,\,\,\,\,\,\,\,\,\,  B=\pi(\mathbb C1\otimes
C^*(\mathbb F_{\infty}))
$$
belong into $B(H)$, with $H$ separable. Let $\{u_n\}$ be the
universal unitaries in $C^*(\mathbb F_{\infty})$. They are clearly a
norm-total sequence. Let
$$
v_n=\pi(u_n\otimes1)\in
A\,\,\,\,\,\,\,\,\,\,\,\,\,\,\,\,\,\,\,\,\,w_n=\pi(1\otimes u_n)\in
B
$$
Now, let $M=A''\in vN(H)$. By hypothesis, there exists a sequence
$\{F_n\}\subseteq\Im_{I_{fin}}$ such that $F_n\rightarrow M$. So
$A\subseteq A''=liminf F_n$. Thus we have
$$
A\subseteq liminf F_n
$$
Moreover
$$
B\subseteq A'=M'=(limsupF_n)'=liminf F_n'
$$
by the commutant theorem (see \cite{Ha-Wi1} th. 3.5). Now we observe
that
$$
\{v_n\}\subseteq U(A)\subseteq U(liminf
F_m)=\underline{lim}Ball(F_m)\cap U(B(H))
$$
where the equality follows from \cite{Ha-Wi1} th.2.6. Let $w(x)$ and
$s^*(x)$ respectively the families of weakly and strong* open
neighborhoods of an element $x\in B(H)$. We have just proved that
for every $n\in\mathbb N$ and $W\in w(v_n)$, one has $W\cap
Ball(F_n)\cap U(B(H))\neq\emptyset$ eventually. Now let $S\in
s^*(u_n)$. By \cite{Ha-Wi1} Lemma 2.4, there exists $W\in w(v_n)$
such that $W\cap Ball(F_m)\cap U(B(H))\subseteq S\cap Ball(F_m)\cap
U(B(H))$. Now, since the first set must be eventually non empty,
also the second one must be the same. This means that we can
approximate (in the strong* topology) $v_n$ with elements in
$U(F_m)$. So let $\{v_{i,n}\}_i\subseteq F_n$ such that
$v_{i,n}\rightarrow^{s^*}v_i$. In a similar way we can find
unitaries $w_{i,n}$ in $F_n'$ such that
$w_{i,n}\rightarrow^{s^*}w_i$. Now let $n$ be fixed, $\pi_{1,n}$ a
representation of $C^*(\mathbb F_{\infty})$ which maps $u_i$ in
$v_{i,n}$ and $\pi_{2,n}$ a representation of $C^*(\mathbb
F_{\infty})$ which maps $u_i$ in $w_{1,n}$. We can find these
representation because the $u_n's$ have no relations among
themselves and because any representation of $G$ extends to a
representation of $C^*(G)$. Notice now that the ranges of these
representations commute, since $v_n\in A$ and $w_n\in B$, and $A,B$
commute. Moreover, the image of $\pi_{1,n}$ belongs into $C^*(F_n)$
and the image of $\pi_{2,n}$ belongs into $C^*(F_n')$. So, by the
universal property in Prop. \ref{universal}, there are unique
representations $\pi_n$ of $C^*(\mathbb
F_{\infty})\otimes_{max}C^*(\mathbb F_{\infty})$ such that
$$
\pi_n(u_i\otimes1)=v_{i,n}\,\,\,\,\,\,\,\,\,and\,\,\,\,\,\,\,\,\,
\pi_n(1\otimes u_i)=w_{i,n}\,\,\,\,\,\,\,\,\,i,n\in\mathbb N
$$
whose image is into $C^*(F_n,F_n')$. Now, since $F_n$ are finite
type $I$ factors, one has $C^*(F_n,F_n')=F_n\otimes_{min}F_n'$ and
thus $\pi_n$ splits: $\pi_n=\sigma_n\otimes\rho_n$, for some
$\sigma_n,\rho_n$ representation of $C^*(\mathbb F_{\infty})$ in
$C^*(F_n,F_n')$. Consequently $||\pi_n(x)||\leq||x||_{min}$ for all
$n\in\mathbb N$ and $x\in C^*(\mathbb F_{\infty})\otimes C^*(\mathbb
F_{\infty})$. On the other hand the sequence $\{\pi_n\}$ converges
to $\pi$ in a strong* pointwise sense (because $\{u_n\}$ is total).
Therefore
$$
||\pi(x)||\leq
liminf||\pi_n(x)||\leq||x||_{min}\,\,\,\,\,\,\,\,\,\forall x\in
C^*(\mathbb F_{\infty})\otimes C^*(\mathbb F_{\infty})
$$
Since $\pi$ is arbitrary, it follows that
$||x||_{max}\leq||x||_{min}$ and the proof of the first implication is complete.\\\\
Notice that we had to work with the strong* topology in order to use
the inequality $||\pi(x)||\leq liminf||\pi_n(x)||$ which fails in
case of weak convergence.\\\\
In order to prove 2. we need two preliminary results
\begin{lem}{\bf (Haagerup-Winslow, \cite{Ha-Wi2} Lemma 4.3)}\label{haag}
Let $A$ be a unital $C^*$-algebra and $\lambda,\rho$ representation
of $A$ in $B(H)$. Assume $\rho$ is faithful and satisfies
$\rho\sim\rho\oplus\rho\oplus...$. Then there exists a sequence
$\{u_n\}\subseteq U(B(H))$ such that
$$
u_n\rho(x)u_n^*\rightarrow^{s^*}\lambda
(x)\,\,\,\,\,\,\,\,\,\,\,\,\,\forall x\in A
$$
\end{lem}
\begin{teo}{\bf (Choi, \cite{Ch} th.7)}\label{ch}
Let $\mathbb F_2$ be the free group with two generators. Then
$C^*(\mathbb F_2)$ has a separating family of finite dimensional
representations.
\end{teo}
(Proof of 2.)\\ By using Choi's theorem and the classical embedding
of $\mathbb F_{\infty}$ into $\mathbb F_2$, we can find a sequence
$\sigma_n$ of finite dimensional representations of $C^*(\mathbb
F_{\infty})$ such that $\sigma=\sigma_1\oplus\sigma_2\oplus...$ is
faithful. Replacing $\sigma$ with the direct sum infinitely many
times of itself, we may assume that
$\sigma\sim\sigma\oplus\sigma\oplus...$. Moreover, by \cite{Ta1}
IV.4.9, $\rho=\sigma\otimes\sigma$ is a faithful representation of
$C^*(\mathbb F_{\infty})\otimes_{min}C^*(\mathbb F_{\infty})$
(because $\rho$ is factorizable). This representation still
satisfies $\rho\sim\rho\oplus\rho\oplus...$. Furthermore, since it
is direct sum of finite dimensional representations, it is separable
and thus we may assume that its image is into $B(H)$, with
$H$ separable.\\
Now, given $M\in vN(H)$, let $\{v_n\},\{w_n\}$ be strong* dense
sequences of unitaries in $Ball(M)$ and $Ball(M')$, respectively.
Let $\{z_n\}$ be the universal unitaries representing $\mathbb
F_{\infty}$ in $C^*(\mathbb F_{\infty})$. Now, by hypothesis and by
using Prop.\ref{universal}, we have a unique representation
$\lambda$ of $C^*(\mathbb F_{\infty})\otimes_{min}C^*(\mathbb
F_{\infty})$ such that
$$
\lambda(z_n\otimes1)=v_n\,\,\,\,\,\,\,\,\,\,\,\,and\,\,\,\,\,\,\,\,\,\,\,\,\lambda(1\otimes
z_n)=w_n\,\,\,\,\,\,\,\,\,\,\,\,\forall n\in\mathbb N
$$
Let us now observe that $M=\lambda(C^*(\mathbb
F_{\infty})\otimes\mathbb C1)''$, since $\{v_n\}$ is dense in
$Ball(M)$. Now, by lemma \ref{haag}, we have unitaries $u_n\in
U(B(H))$ such that
$$
u_n\rho(x)u_n^*\rightarrow^{s^*}\lambda(x)\,\,\,\,\,\,\,\,\,\,\,\,\,\,\,\,\,\,\,\,\,\,\,\,\forall
x\in C^*(\mathbb F_{\infty})\otimes_{min}C^*(\mathbb F_{\infty})
$$
Define
$$
M_n=u_n\rho(C^*(\mathbb F_{\infty})\otimes\mathbb C1)''u_n^*
$$
then
$$
u_n\rho(\mathbb C1\otimes C^*(\mathbb F_{\infty}))u_n^*\subseteq
M_n'
$$
So we have (by using Rem.\ref{strong})
$$
\lambda(C^*(\mathbb F_{\infty})\otimes\mathbb
C1)=liminfu_n\rho(C^*(\mathbb F_{\infty})\otimes\mathbb
C1)u_n^*\subseteq liminf M_n
$$
In a similar way, we obtain
$$
\lambda(\mathbb C1\otimes C^*(\mathbb F_{\infty}))\subseteq liminf
M_n'
$$
Now, by \cite{Ha-Wi1} th.2.3, $liminf M_a$ is always a von Neumann
algebra, and thus the previous inclusions hold by passing to the
strong closure:
$$
M=\lambda(C^*(\mathbb F_{\infty})\otimes\mathbb C1)''\subseteq
liminf M_n
$$
and
$$
M'=\lambda(\mathbb C1\otimes C^*(\mathbb F_{\infty}))''\subseteq
liminf M_n'
$$
Now, applying the commutant theorem $(liminf M_a)'=limsup M_a'$ (see
\cite{Ha-Wi1}, th.3.5), we have $M_n\rightarrow M$. Now, we observe
that $\rho$ is a type $I$ representation, since it is direct sum of
finite dimensional representations, and thus $M_n\in vN_I(H)$. Thus
we have just proved that $vN_I(H)$ is dense in $vN(H)$. In
particular $vN_{inj}(H)$ is dense in $vN(H)$. Now it has been
already proved by Haagerup and Winslow that $vN_{inj}$ and $\Im(H)$
(factors into $B(H))$ are $G_{\delta}$ and $\Im(H)$ is dense. On the
other hand, $vN(H)$ is a Polish space and hence a Baire's space. So,
also the intersection $vN_{inj}(H)\cap\Im(H)=\Im_{inj}(H)$ must be
dense.
\end{proof}
\end{teo}
Notice that in the proof of 2. we have used the hypothesis only to
apply Prop.\ref{universal}. We need it to have $\lambda$ and $\rho$
defined on the same $C^*$-algebra and so apply Lemma \ref{haag}.
\begin{nota}
One can ask what groups $G$ satisfy \emph{Kirchberg's property}
$$
C^*(G)\otimes_{min}C^*(G)=C^*(G)\otimes_{max}C^*(G)
$$
or the \emph{reduced Kirchberg's property}
$$
C_r^*(G)\otimes_{min}C_r^*(G)=C_r^*(G)\otimes_{max}C_r^*(G)
$$
where $C_r^*(G)$ is the reduced $C^*$-algebra of $G$, i.e. the
C*-algebra generated by the image of the left regular representation
on $l^2(G)$. Let us denote by $K$ and $K_r$ respectively the classes
of group which satisfy \emph{Kirchberg's property} and \emph{reduced
Kirchberg's property}. It follows from a more general result by
Conti and Hamhalter (see \cite{Co-Ha}) that
$$
K_r\cap\{i.c.c.\,\,\,groups\}=\{amenable\,\,\,groups\}
$$
What can we say about $K$? Are there any non-amenable examples?
\end{nota}
\section{Lance's WEP and QWEP conjecture}
Kirchberg's theorem \ref{kirchi} shows an interesting and unexpected
link between von Neumann algebras and $C^*$-algebras. Kirghberg
himself, in \cite{Ki} again, has found another interesting link
between them; more precisely: Connes' embedding conjecture is a
particular case of a conjecture regarding the structure of
$C^*$-algebras: QWEP conjecture. We remind the reader that a
$C^*$-algebra is QWEP if it is a quozient of a $C^*$-algebra with
Lance's WEP. So it is natural to ask if there is a direct relation
between Connes' embedding conjecture and WEP. N.P. Brown has found
in \cite{Br} that Connes' embedding conjecture is equivalent
to the analogue of Lance's WEP for separable type $II_1$ factors.\\
In order to give some details about this let us firstly recall the
following
\begin{defin}
Let $A,B$ be $C^*$-algebras and $\phi:A\rightarrow B$ a linear map.
For every $n\in\mathbb N$ we can define a map
$\phi_n:M_n(A)\rightarrow M_n(B)$ by setting
$$
\phi_n[a_{ij}]=[\phi(a_{ij})]
$$
$\phi$ is called completely positive if $\phi_n$ is positive for
every $n$.
\end{defin}
\begin{nota}
Any *homomorphism between two $C^*$-algebras is automatically c.p.
Indeed it is clearly positive. On the other hand $\phi_n$ can be
described as $\phi\otimes Id_n$ and thus it is still a
*homomorphism, since tensor product of *homomorphisms is still a
*homomorphism.
\end{nota}
\begin{defin}
Let $A\subseteq B$ be two $C^*$-algebras. We say that $A$ is weakly
cp complemented in $B$ if there exists a unital completely positive
map $\phi:B\rightarrow A^{**}$ such that $\phi|_A=Id_A$.
\end{defin}
\begin{defin}
We say that a $C^*$-algebra $A$ has the WEP (weak expectation
property) if it is weakly cp complemented in $B(H)$ for a faithful
representation $A\subseteq B(H)$.
\end{defin}
This property is been introduced by Lance in \cite{La}, where he
proved also that this definition does not depend on the choice of
the faithful representation of $A$.
\begin{defin}
We say that $A$ has QWEP if it is a quotient of a $C^*$-algebra with
WEP.
\end{defin}
Here is QWEP conjecture, regarding the structure of a $C^*$-algebra.
\begin{con}{\bf (QWEP conjecture)}
Every $C^*$-algebra is QWEP.
\end{con}
The unexpected theorem by Kirghberg is
\begin{teo}{\bf (Kirchberg)}
The following statements are equivalent
\begin{enumerate}
\item Connes' embedding conjecture is true.
\item QWEP conjecture is true for separable von Neumann algebras.
\end{enumerate}
\end{teo}
A proof of this theorem can be found in the original paper by
Kirchberg \cite{Ki} or also in \cite{Oz}. Now we prefer to focus on
an easier and equally interesting topic: the von Neumann algebraic
analogue of Lance's WEP and the proof of Brown's theorem. What
follows is just a rewriting of Brown's paper \cite{Br}.
\begin{defin}
Let $M\subseteq B(H)$ a von Neumann algebra and $A\subseteq M$ a
weakly dense $C^*$-subalgebra. We say that $M$ has a weak
expectation relative to $A$ if there exists a u.c.p. map
$\Phi:B(H)\rightarrow M$ such that $\Phi(a)=a$, for all $a\in A$.
\end{defin}
\begin{nota}
The notion of injectivity for von Neumann algebras can be given also
in the following way: $M\subseteq B(H)$ is injective if there exists
a u.c.p. map $\Phi:B(H)\rightarrow M$ such that $\Phi(x)=x$, for all
$x\in M$. So weak expectation relative property is something less
than injectivity. Actually something more precise holds: Brown's
theorem can be read by saying that weak expectation relative
property is the \emph{limit property of injectivity}. We can clarify
this interpretation after enunciating the following
\end{nota}
\begin{teo}{\bf (Brown, \cite{Br} Th.1.2)}\label{wep}
For a separable type $II_1$ factor $M$ the following conditions are
equivalent:
\begin{enumerate}
\item $M$ is embeddable into $R^{\omega}$.
\item $M$ has a weak expectation relative to some weakly dense
subalgebra.
\end{enumerate}
\end{teo}
We can now clarify the interpretation of the weak expectation
relative property as limit of injectivity.
\begin{cor}
For a separable type $II_1$ factor the following conditions are
equivalent:
\begin{enumerate}
\item $M$ has a weak expectation relative property.
\item $M$ is Effros-Marechal limit of injective factors.
\end{enumerate}
\begin{proof}
It is an obvious consequence of Th.\ref{wep} and Th.\ref{haa}.
\end{proof}
\end{cor}
Our purpose is to present the original proof of Th.\ref{wep}. We need some preliminary result.\\\\
For the rest of the chapter let $A$ be a separable $C^*$-algebra.
This hypothesis is not necessary, but it is convenient.
\begin{defin}
A tracial state on $A$ is map $\tau:A_+\rightarrow[0,\infty]$ such
that
\begin{enumerate}
\item $\tau(x+y)=\tau(x)+\tau(y)$, for all $x,y\in A_+$
\item $\tau(\lambda x)=\lambda\tau(x)$, for all $\lambda\geq0,x\in
A_+$
\item $\tau(x^*x)=\tau(xx^*)$ for all $x\in A$
\item $\tau(1)=1$
\end{enumerate}
It clearly extends to a positive functional on the whole $A$.
\end{defin}
\begin{defin}
A tracial state $\tau$ on $A\subseteq B(H)$ is called
\emph{invariant mean} if there exists a state $\psi$ on $B(H)$ such
that
\begin{enumerate}
\item $\psi(uTu^*)=\psi(T)$, for all $u\in U(A)$ and $T\in B(H)$
\item $\psi|_A=\tau$
\end{enumerate}
\end{defin}
\begin{nota}
A consequence of Th.\ref{invariant} is that the notion of invariant
mean does not depend on the choice of the faithful representation
$A\subseteq B(H)$.
\end{nota}
In order to prove Brown's theorem we need a characterization of
invariant means.\\\\
We recall the following well-known
\begin{teo}{\bf (Powers-St{\o}rmer inequality, \cite{Po-St})}
Let $h,k\in L^1(B(H))_+$. Then
$$
||h-k||_2^2\leq||h^2-k^2||_1
$$
where $||\cdot||_i$ stands for the $L^i$ norm on $L^1(B(H))$ with
respect to the canonical unbounded trace Tr. In particular, if $u\in
U(B(H))$ and $h\geq0$ has finite rank, then
$$
||uh^{1/2}-h^{1/2}u||_2=||uh^{1/2}u^*-h^{1/2}||_2\leq||uhu^*-h||_1^{1/2}
$$
\end{teo}
\begin{lem}\label{brown1}
Let $H$ be a separable Hilbert space and $h\in B(H)$ a positive,
finite rank operator with rational eigenvalues and $Tr(h)=1$. Then
there exists a u.c.p. map $\Phi:B(H)\rightarrow M_q(\mathbb C)$ such
that
\begin{enumerate}
\item $tr(\Phi(T))=Tr(hT)$, for all $T\in B(H)$
\item $|tr(\Phi(uu^*)-\Phi(u)\Phi(u^*))|<2||uhu^*-h||_1^{1/2}$, for
all $u\in U(B(H))$
\end{enumerate}
Here $tr$ stands for the normalized trace on $M_q(\mathbb C)$.
\begin{proof}
Let $v_1,... v_k\in H$ be the eigenvectors of $H$ and
$\frac{p_1}{q},...\frac{p_k}{q}$ the corresponding eigenvalues. Thus
\begin{enumerate}
\item $hv_i=\frac{p_i}{q}$
\item $\sum_{i=1}^k\frac{p_i}{q}=tr(h)=1$. It follows that $\sum
p_i=q$
\end{enumerate}
Let $\{w_m\}$ be any orthonormal basis for $H$. Consider the
orthogonal subset of $H\otimes H$:
$$
B=\{v_1\otimes w_1,...v_1\otimes v_{p_1}\}\cup\{v_2\otimes w_1,...
v_2\otimes w_{p_2}\}\cup...\cup\{v_k\otimes w_1,... v_k\otimes
w_{p_k}\}
$$
Let $V$ be the subspace of $H\otimes H$ spanned by $B$ and
$P:H\otimes H\rightarrow V$ the orthogonal projection. Let $T\in
B(H)$, the following formula holds
$$
Tr(P(T\otimes1)P)=\sum_{i=1}^kp_i<Tv_i,v_i>
$$
Indeed $P(T\otimes1)P$ is representable (in the basis $B$) by a
$q\times q$ block diagonal matrix whose blocks have dimension $p_i$
with entries $ETE$, where $E:H\rightarrow span\{v_1,...v_k\}$ is the
projection. Now define $\Phi:B(H)\rightarrow M_q(\mathbb C)$ by
setting $\Phi(T)=P(T\otimes1)P$. We have
$$
tr(\Phi(T))=\frac{1}{q}Tr(P(T\otimes1)P)=\sum_{i=1}^k<T\frac{p_i}{q}v_i,v_i>=\sum_{i=1}^k<Thv_i,v_i>=Tr(Th)
$$
Moreover $\Phi$ is u.c.p. So the first assertion is proved.\\
Now, by writing down the matrix of $P(T\otimes1)P(T^*\otimes1)P$ in
the basis $B$ we have
$$
Tr(P(T\otimes1)P(T^*\otimes1)P)=\sum_{i,j=1}^k|T_{i,j}|^2min(p_i,p_j)
$$
where $T_{i,j}=<Tv_j,v_i>$. Analogously, by writing down the
matrices of $h^{1/2}T, Th^{1/2}$ and $h^{1/2}Th^{1/2}T^*$ in any
orthonormal basis which begins with $\{v_1,... v_k\}$ we have
$$
Tr(h^{1/2}Th^{1/2}T^*)=\sum_{i,j=1}^k\frac{1}{q}(p_ip_j)^{1/2}|T_{i,j}|^2
$$
By using these formulas, we can make the following preliminary
calculation
$$
|Tr(h^{1/2}Th^{1/2}T^*)-tr(\Phi(T)\Phi(T^*))|=
$$
$$
=|\sum_{i,j=1}^k\frac{1}{q}(p_ip_j)^{1/2}|T_{i,j}|^2-\frac{1}{q}Tr(P(T\otimes1)P(T^*\otimes1)P)|=
$$
$$
=|\sum_{i,j=1}^k\frac{1}{q}|T_{i,j}|^2((p_ip_j)^{1/2}-min(p_i,p_j)|\leq
$$
by using $min(p_i,p_j)\leq p_i$
$$
\leq\sum_{i,j=1}^k\frac{1}{q}|T_{i,j}|^2p_i^{1/2}|p_j^{1/2}-p_i^{1/2}|\leq
$$
by using the Holder inequality
$$
\leq
(\sum_{i,j=1}^k\frac{1}{q}|T_{i,j}|^2p_i)^{1/2}(\sum_{i,j=1}^k\frac{1}{q}|T_{i,j}|^2(p_i^{1/2}-p_j^{1/2}))^{1/2}=
$$
$$
=||Th^{1/2}||_2||h^{1/2}T-Th^{1/2}||_2=
$$
suppose now that $T\in U(B(H))$, so that
$||Th^{1/2}||_2=||h^{1/2}||_2=1$
$$
=||h^{1/2}T-Th^{1/2}||_2=||Th^{1/2}T^*-h^{1/2}||_2\leq
$$
by using the Powers-St{\o}rmer inequality
$$
\leq||ThT^*-T||_1^{1/2}
$$
Now we can prove the second assertion. Indeed we have
$$
|Tr(\Phi(TT^*)-\Phi(T)\Phi(T^*))|\leq
$$
by using the triangle inequality and the previous calculation
$$
|1-Tr(h^{1/2}Th^{1/2}T^*)|+||ThT^*-h||_1^{1/2}=
$$
$$
=|Tr(ThT^*)-Tr(h^{1/2}Th^{1/2}T^*)|+||ThT^*-h||_1^{1/2}=
$$
$$
=|Tr((Th^{1/2}-h^{1/2}T)h^{1/2}T^*)|+||ThT^*-h||_1^{1/2}\leq
$$
by using the Cauchy-Schwarz inequality
$$
\leq||h^{1/2}T^*||_2||Th^{1/2}-h^{1/2}T||_2+||ThT^*-h||_1^{1/2}
$$
So the assertion follows by using $T\in U(B(H)), Tr(h)=1$ and by
applying the Powers-St{\o}rmer inequality once more.
\end{proof}
\end{lem}
We recall a classical theorem by Choi
\begin{teo}{\bf (Choi, \cite{Ch2})}\label{choi}
Let $A,B$ be two $C^*$-algebras and $\Phi:A\rightarrow B$ a u.c.p.
map. Then
$$
\{a\in A:\Phi(aa^*)=\Phi(a)\Phi(a^*),\Phi(a^*a)=\Phi(a^*)\Phi(a)\}=
$$
$$
=\{a\in A:\Phi(ab)=\Phi(a)\Phi(b),\Phi(ba)=\Phi(b)\Phi(a),\forall
b\in A\}
$$
\end{teo}
Here is the characterization of invariant means. Other ways to
characterize them are in \cite{Br2}, Th.3.1, and in \cite{Oz}, Th.
6.1.
\begin{teo}\label{invariant}
Let $\tau$ be a tracial state on $A\subseteq B(H)$. Then the
followings are equivalent:
\begin{enumerate}
\item $\tau$ is an invariant mean.
\item There exists a sequence of u.c.p. maps $\Phi_n:A\rightarrow
M_{k(n)}(\mathbb C)$ such that
\begin{enumerate}
\item $||\Phi_n(ab)-\Phi_n(a)\Phi_n(b)||_2\rightarrow0$ for all
$a,b\in A$
\item $\tau(a)=lim_{n\rightarrow\infty}tr(\Phi_n(a))$, for all $a\in
A$
\end{enumerate}
\item For any faithful representation $\rho:A\rightarrow B(H)$ there exists
a u.c.p. map $\Phi:B(H)\rightarrow\pi_{\tau}(A)''$ such that
$\Phi(\rho(a))=\pi_{\tau}(a)$, for all $a\in A$, where $\pi_{\tau}$
stands for the GNS representation associated to $\tau$.
\end{enumerate}
\begin{proof}
($1\Rightarrow 2$)\\
Let $\tau$ be an invariant mean with respect to the faithful
representation $\rho:A\rightarrow B(H)$. Thus we can find a state
$\psi$ on $B(H)$ which extends $\tau$ and such that
$\psi(uTu^*)=\psi(u)$, for all $u\in U(A)$ and for all $T\in B(H)$.
Since the normal states are dense in $B(H)$ and they are represented
in the form $Tr(h\cdot)$, with $h\in L^1(B(H))$, we can find a net
$h_{\lambda}\in L^1(B(H))$ such that
$Tr(h_{\lambda}T)\rightarrow\psi(T)$, for all $T\in B(H)$. Moreover
we remind that $h_{\lambda}$ is positive and has trace $1$. Now,
since $\psi(uTu^*)=\psi(u)$, it follows that
$Tr(uh_{\lambda}u^*T)=Tr(h_{\lambda}u^*Tu)\rightarrow\psi(u^*Tu)=\psi(T)$
and thus $Tr(h_{\lambda}T)-Tr((uh_{\lambda}u^*)T)\rightarrow 0$, for
all $T\in B(H)$, i.e. $h_{\lambda}-uh_{\lambda}u^*\rightarrow0$ in
the weak topology on $L^1(B(H))$. Now let $\{U_n\}$ be an increasing
family of finite sets of unitaries whose union have dense linear
span in $A$ and $\varepsilon=\frac{1}{n}$. Let $U_{n}=\{u_1,...
u_n\}$. Fixed $n$, let us consider the convex hull of the set
$\{u_1h_{\lambda}u_1^*-h_{\lambda},...u_nh_{\lambda}u_n^*-h_{\lambda}\}$.
Its weak closure contains $0$ (because of the previous observation)
and coincide with the $1$-norm closure, by the Hahn-Banach
separation theorem. Thus there exists a convex combination of
$h_{\lambda}$'s, say $h$, such that
\begin{enumerate}
\item $Tr(h)=1$
\item $||uhu^*-h||_1<\varepsilon, \forall u\in U_n$
\item $|Tr(uh)-\tau(u)|<\varepsilon, \forall u\in U_n$
\end{enumerate}
Moreover, since finite rank operators are norm dense in $L^1(B(H))$,
we can suppose that $h$ is finite rank with rational eigenvalues.
Now we can apply Lemma \ref{brown1} in order to construct a sequence
of u.c.p. maps $\Phi_n:B(H)\rightarrow M_{k(n)}(\mathbb C)$ such
that
\begin{enumerate}
\item $Tr(\Phi_n(u))\rightarrow\tau(u)$
\item $|Tr(\Phi_n(uu^*))-\Phi_n(u)\Phi_n(u^*)|\rightarrow0$
\end{enumerate}
for every unitary in a countable set whose linear span is dense in
$A$. So we have obtained the thesis for unitaries. The second
property holds for any $a\in A$, by passing to linear combinations.
In order to obtain the first one, we observe that
$\Phi_n(uu^*)-\Phi_n(u)\Phi_n(u^*)\geq0$ and thus the following
inequality holds
$$
||1-\Phi_n(u)\Phi_n(u^*)||_2^2\leq||1-\Phi_n(u)\Phi_n(u^*)||tr(\Phi_n(uu^*)-\Phi_n(u)\Phi_n(u^*))
$$
and the right hand side tends to zero. Now define
$\Phi=\oplus\Phi_n:A\rightarrow\Pi M_{k(n)}(\mathbb C)\subseteq
l^{\infty}(R)$ and compose with the quotient map
$p:l^{\infty}(R)\rightarrow R^{\omega}$. The previous inequality
shows that if $u$ is a unitary such that
$||\Phi_n(uu^*)-\Phi_n(u)\Phi_n(u^*)||_2\rightarrow0$ and
$||\Phi_n(u^*u)-\Phi_n(u^*)\Phi_n(u)||_2\rightarrow0$, then $u$
falls in the multiplicative domain of $p\circ\Phi$. But such
unitaries have dense linear span in $A$ and hence the whole $A$
falls in the multiplicative domain of $p\circ\Phi$ (by Choi's
theorem \ref{choi}). By definition of ultraproduct this just means
that $||\Phi_n(ab)-\Phi_n(a)\Phi_n(b)||_2\rightarrow0$, for all
$a\in A$.
\\\\
($2\Rightarrow3$)\\
Let $\Phi_n:A\rightarrow M_{k(n)}(\mathbb C)$ be a sequence of
u.c.p. maps with the properties stated in the theorem. Identify each
$M_{k(n)}(\mathbb C)$ with a unital subfactor of $R$ and we can
define a u.c.p map $\tilde{\Phi}:A\rightarrow l^{\infty}(R)$ by
$x\rightarrow(\Phi_n(x))_n$. Since the $\Phi_n's$ are asymptotically
multiplicative in 2-norm one get a $\tau$-preserving *homomorphism
$A\rightarrow R^{\omega}$ by composing with the quotient map
$p:l^{\infty}(R)\rightarrow R^{\omega}$. Note that the weak closure
of $p\circ\tilde{\Phi}(A)$ into $R^{\omega}$ is isomorphic to
$\pi_{\tau}(A)''$. Thus we are in the following situation
$$
\xymatrix{A\ar[r]^{\tilde{\Phi}}\ar[d]^{\rho} & l^{\infty}(R)\ar[r]^p & R^{\omega} & \supseteq & \overline{p\circ\tilde{\Phi}(A)}^w\cong \pi_{\tau}(A)''\\
B(H)\ar[d]^i\\
B(K)}
$$
where $K$ is a representing Hilbert space for $l^{\infty}(R)$ and
$i$ is a natural embedding ($K$ cannot be separable). Now
$l^{\infty}(R)$ is injective and let $E:B(K)\rightarrow
l^{\infty}(R)$ a surjective projection of norm 1. Moreover let
$F:R^{\omega}\rightarrow\pi_{\tau}(A)''$ a conditional expectation
(see \cite{Ta1}, Prop.2.36). Thus we are in the following situation
$$
\xymatrix{A\ar[r]^{\tilde{\Phi}}\ar[d]^{\rho} & l^{\infty}(R)\ar[r]^p & R^{\omega}\ar[r]^F & \pi_{\tau}(A)''\cong\overline{p\circ\tilde{\Phi}(A)}^w\\
B(H)\ar[d]^i\\
B(K)\ar[uur]^E}
$$
Define $\Phi:B(H)\rightarrow\pi_{\tau}(A)''$ by setting $\Phi=FpEi$.
Clearly $\Phi(\rho(a))=\pi_{\tau}(a)$.\\\\
($3\Rightarrow1$)\\
The hypothesis $\Phi(a)=\pi_{\tau}(a)$ guarantees that $\Phi$ is
multiplicative on $A$. By Choi's theorem \ref{choi} it follows that
$\Phi(aTb)=\pi_{\tau}(a)\Phi(T)\pi_{\tau}(b)$, for all $a,b\in A,
T\in B(H)$. Let $\tau''$ be the vector trace on $\pi_{\tau}(A)''$
and consider $\tau''\circ\Phi$. Clearly it extends $\tau$. Moreover
it is invariant under the action of $U(A)$, indeed
$$
(\tau''\circ\Phi)(u^*Tu)=\tau''(\pi_{\tau}(u)^*\Phi(T)\pi_{\tau}(u))=\tau''(\Phi(T))=(\tau''\circ\Phi)(T)
$$
Hence $\tau$ is an invariant mean.
\end{proof}
\end{teo}
Another preliminary but very nice result is the following
\begin{prop}\label{densealgebra}
Let $M$ be a separable type $II_1$ factor. There exists a
*-monomorphism $\rho:C^*(\mathbb F_{\infty})\rightarrow M$ such that
$\rho(C^*(\mathbb F_{\infty}))$ is weakly dense in $M$.
\begin{proof}
We first observe that $C^*(\mathbb F_{\infty})$ is inductive limit
of free products of itself. It can be imagined by partitioning the
set of generators in a sequence of countable set (one can do it
because $|\mathbb N|=|\mathbb N\times\mathbb N|$). Let $\{X_n\}$
such a sequence. Define $A_n=C^*(X_1,...,X_n)$. Clearly one has
$A_n=A_{n-1}*C^*(X_n)$, where $*$ stands for the free product with
amalgamation over the scalar. Moreover $C^*(X_n)\cong C^*(\mathbb
F_{\infty})$, and then $A_n\cong A_{n-1}*C^*(\mathbb F_{\infty})$.
Now let $A$ be the inductive limit of the $A_n's$. Clearly
$A=\bigcup A_n=C^*(X_1,X_2,...)\cong C^*(F_{\infty})$. Now, by
Choi's theorem \ref{ch} we can find a sequence of integers
$\{k(n)\}$ and a unital
*-monomorphism $\sigma:A\rightarrow\Pi M_{k(n)}(\mathbb C)$. Note
that we may naturally identify each $A_i$ with a subalgebra of $A$
and hence, restricting $\sigma$ to this copy, get an injection of
$A_i$ into $\Pi M_{k(n)}(\mathbb C)$. Now we can prove the existence
of a sequence of unital *homomorphism $\rho_i:A_i\rightarrow M$ such
that
\begin{enumerate}
\item Each $\rho_i$ is injective;
\item $\rho_{i+1}|_{A_i}=\rho_i$ where we identify $A_i$ with the
''left side'' of $A_i*C^*(\mathbb F_{\infty})=A_{i+1}$;
\item The union of $\{\rho_i(A_i)\}$ is weakly dense in $M$.
\end{enumerate}
After finding the $\rho_i's$, it will be enough to define $\rho$ as
union of those ones.\\
We first choose an increasing sequence of projections of $M$ such
that $\tau_M(p_i)\rightarrow1$. Then we define the orthogonal
projections $q_n=p_n-p_{n-1}$ and consider the type $II_1$ factors
$Q_i=q_iMq_i$. Now, by the division property of type $II_1$ factors,
we can find a unital embedding $\Pi M_{k(n)}\rightarrow Q_i\subseteq
M$. By composing with $\sigma$, we get a sequence of embeddings
$A\rightarrow M$, which will be denoted by $\sigma_i$. Now $p_iMp_i$
is weakly separable and thus there is a countable total family of
unitaries. Hence we can find a *homomorphism $\pi_i:C^*(\mathbb
F_{\infty})\rightarrow p_iMp_i$ with weakly dense range (take the
generators of $\mathbb F_{\infty}$ into $C^*(\mathbb F_{\infty})$
and map them into that total family of unitaries). Now we define
$$
\rho_1=\pi_1\oplus(\bigoplus_{j\geq2}\sigma_j|_{A_1}):A_1\rightarrow
p_1Mp_1\oplus(\Pi_{j\geq2}Q_j)\subseteq M
$$
It is a *monomorphism, since each $\sigma_i$ is already faithful on
the whole $A$. Now define a *homomorphism
$\theta_2:A_2=A_1*C^*(\mathbb F_{\infty})\rightarrow p_2Mp_2$ as the
free product of the *homomorphism $A_1\rightarrow p_2Mp_2$,
$x\rightarrow p_2\rho_1(x)p_2$, and $\pi_2:C^*(\mathbb
F_{\infty})\rightarrow p_2Mp_2$. We then put
$$
\rho_1=\theta_2\oplus(\bigoplus_{j\geq3}\sigma_j|_{A_2}):A_2\rightarrow
p_2Mp_2\oplus(\Pi_{j\geq3}Q_j)\subseteq M
$$
Clearly $\rho_2|_{A_1}=\rho_1$. In general, we construct a map
$\theta_{n+1}:A_n*C^*(\mathbb F_{\infty})\rightarrow
p_{n+1}Mp_{n+1}$ as the free product of the cutdown (by $p_{n+1}$)
of $\rho_n$ and $\pi_n$. This map need not be injective and hence we
take a direct sum with $\oplus_{j\geq n+2}\sigma_j|_{A_{n+1}}$ to
remedy this deficiency. These maps have all the required properties
and hence the proof is complete (note that the last property follows
from the fact that the range of each $\theta_n$ is weakly dense in
$p_{n+1}Mp_{n+1}$).
\end{proof}
\end{prop}
Now we can prove Brown's theorem
\begin{teo}{\bf (Brown)}\label{marrone}
Let $M$ be a separable type $II_1$ factor. The followings are
equivalent:
\begin{enumerate}
\item $M$ is embeddable into $R^{\omega}$.
\item $M$ has the weak expectation property relative to some weakly
dense subalgebra.
\end{enumerate}
\begin{proof}
($1\Rightarrow2$)\\
Let $M$ be embeddable into $R^{\omega}$. By Prop.\ref{densealgebra},
we may identify $C^*(\mathbb F_{\infty})$ with a weakly dense
subalgebra $A$ of $M$. We want to prove that $M$ has the weak
expectation property relative to $A$. Let $\tau$ the unique
normalized trace on $M$, more precisely we will prove that
$\pi_{\tau}(M)$ has the weak expectation property relative to
$\pi_{\tau}(A)$. Indeed $\tau$ is faithful and w-continuous and
hence $\pi_\tau(M)$ and $\pi_\tau(A)$ are respectively copies of $M$
and $A$ and $\pi_\tau(A)$ is still weak dense in $\pi_\tau(M)$. We
first prove that $\tau|_A$ is an invariant mean. Take $\{u_n\}$
universal generators of $\mathbb F_{\infty}$ into $A$. Let $n$ be
fixed, since $u_n\in R^\omega$ it is $||\cdot||_2$-limit of
unitaries in $R$; on the other hand, the unitary matrices are weakly
dense in $U(R)$ and hence they are $||\cdot||_2$-dense in $U(R)$
(since w-closed convex subsets coincide with the
$||\cdot||_2$-closed convex ones (see \cite{Jo})). Thus we can find
a sequence of unitary matrices which converges to $u_n$ in norm
$||\cdot||_2$. Let $\sigma$ be the mapping which sends each $u_n$ to
such a sequence. Since the $u_n$'s have no relations, we can extend
$\sigma$ to a
*homomorphism $\sigma:C^*(\mathbb F_{\infty})\rightarrow\Pi
M_k(\mathbb C)\subseteq l^{\infty}(R)$. Let
$p:l^{\infty}(R)\rightarrow R^{\omega}$ be the quotient mapping. By
the $2$-norm convergence we have $(p\circ\sigma)(x)=x$ for all
$x\in\mathbb C^*(\mathbb F_{\infty})$. Let $p_n:\Pi M_k(\mathbb
C)\rightarrow M_n(\mathbb C)$ be the projection, by the definition
of trace in $R^{\omega}$, we have
$$
\tau(x)=lim_{n\rightarrow\omega}tr_n(p_n(\sigma(x)))
$$
where $tr_n$ is the normalized trace on $M_n(\mathbb C)$. Now we can
apply \ref{invariant},2) by setting $\phi_n=p_n\circ\sigma$ (they
are u.c.p. since they are *homomorphisms) and conclude that
$\tau|_A$ is an invariant mean. Now consider $\pi_{\tau}(M)\subseteq
B(H)$ and $\pi_{\tau}(A)=\pi_{\tau|_A}(A)\subseteq B(H)$. By
Th.\ref{invariant} there exists a u.c.p. map
$\Phi:B(H)\rightarrow\pi_{\tau}(A)''=\pi_{\tau}(M)$ such that
$\Phi(a)=\pi_{\tau}(a)$. Thus $M$ has the weak expectation
property relative to $C^*(\mathbb F_{\infty})$.\\\\
($2\Rightarrow1$)\\
Let $A\subseteq M\subseteq B(H)$ be weakly dense and
$\Phi:B(H)\rightarrow M$ a u.c.p. map which fixes $A$. Let $\tau$ be
the unique normalized trace on $M$. After identifying $A$ with
$\pi_\tau(A)$, we are under the hypothesis of \ref{invariant}.3) and
thus $\tau|_A$ is an invariant mean. By Th.\ref{invariant} it
follows that there exists a sequence $\phi_n:A\rightarrow
M_{k(n)}(\mathbb C)$ such that
\begin{enumerate}
\item $||\phi_n(ab)-\phi_n(a)\phi_n(b)||_2\rightarrow0$ for all
$a,b\in A$
\item $\tau(a)=lim_{n\rightarrow\infty}tr_n(\phi_n(a))$, for all
$a\in A$
\end{enumerate}
Let $p:l^{\infty}(R)\rightarrow R^{\omega}$ be the quotient mapping.
The previous properties guarantee that the u.c.p. mapping
$A\rightarrow R^{\omega}$, $\Phi:x\rightarrow p(\{\phi_n(x)\})$ is a
*homomorphism which preserves $\tau|_A$. It follows it mapping is injective too, since
$\Phi(x)=0\Rightarrow\Phi(x^*x)=0\Rightarrow\tau(x^*x)=0\Rightarrow
x^*x=0\Rightarrow x=0$. Observe now that the weak closure of $A$
into $R^{\omega}$ is isomorphic to $M$ (they are algebraically
isomorphic and have the same trace) and hence $M$ embeds into
$R^{\omega}$.
\end{proof}
\end{teo}
\section{A few words about other approaches}
\subsection{Relation with Hilbert's 17th problem}
The original version of Hilbert's
17th problem is very simple. Let us recall that $\mathbb
R[x_1,...,x_n]$ denotes the ring of polynomials with $n$
indeterminates and real coefficients and $\mathbb R(x_1,...,x_n)$
denotes the quotient field of $\mathbb R[x_1,...,x_n]$.
\begin{prob}{\bf (Hilbert's 17th)}
Given a polynomial $f\in\mathbb R[x_1,...,x_n]$ which is
non-negative for all substitutions $(x_1,...,x_n)\in\mathbb R^n$. Is
it possible to write $f$ as sum of squares of elements in $\mathbb
R(x_1,...,x_n)$?
\end{prob}
The affirmative answer was given by Emil Artin, in 1927 (see
\cite{Ar}). He gave a very abstract solution. Actually, now we have
also an algorithm to construct such a decomposition. It has been
recently found by Delzell (see
\cite{De}).\\
More recently, many mathematicians have looked for generalizations
of the problem. The first and most intuitive one is the following
\begin{prob}
Are the matrices with entries in $\mathbb R[x_1,...,x_n]$ which are
always positive semidefinite (i.e. for all substitutions
$(x_1,...,x_n)\in\mathbb R^n$) sum of squares of symmetric matrices
with entries in $\mathbb R(x_1,...,x_n)$?
\end{prob}
Also in this case an affirmative answer was given independently by
Gondard-Ribenoim (see \cite{Go-Ri}) and Procesi-Schacher (see
\cite{Pr-Sc}). Also in this case, for a constructive solution one
had to wait for thirty years: it has been just found by Hillar and
Nie in
2006 (see \cite{Hi-Ni}).\\\\
Other generalizations come from Geometry and Operator Algebras. Let
us recall the following
\begin{defin}
An $n$-manifold $M$ is called irreducible if for any embedding of
$S^{n-1}$ into $M$ there exists an embedding of $B^n$ into $M$ such
that the image of the boundary $\partial B^n$ coincides with the
image of $S^{n-1}$.
\end{defin}
\begin{prob}{\bf (Geometric
version)}\label{hilbert3} Let $M$ be a paracompact irreducible
analytic manifold and $f:M\rightarrow\mathbb R$ a non-negative
analytic function. Can $f$ be written as a sum of squares of
meromorphic functions?
\end{prob}
We recall that meromorphic functions are analytic functions in the
whole domain except a set of isolates points which are poles. So,
rational functions are meromorphic and one can recognize a
generalization of Hilbert's 17th problem.\\
This problem was solved by Ruiz (see \cite{Ru}) in the case of
compact manifold. In the generale case there are lots of approaches
in course, but a complete solution is known only for $n=2$
(see \cite{Ca}).\\\\
Now we want to describe briefly the formulation of the problem in
terms of Operator Algebras. It is due to R\u adulescu, who proved in
\cite{Ra2} the
equivalence between it and Connes' embedding conjecture.\\
The basic idea is to generalize analytic functions with formal
series. Let $Y_1,...,Y_n$ be $n$ inderminates. We set
$$
\mathcal I_n=\{(i_1,...,i_p\},p\in\mathbb N,
i_1,...,i_p\in\{1,...,n\}\}
$$
If $I=\{i_1,...,i_p\}\in\mathcal I_n$, we set
$Y_I=Y_{i_1}\cdot...\cdot Y_{i_p}$. Let
$$
V=\{\sum_{I\in\mathcal I_n}a_IY_I, a_I\in\mathbb C, ||\sum
a_IY_I||_R=\sum|a_I|R^{|I|}<\infty, \forall R>0\}
$$
So, for any $R>0$, we have a norm on $V$. R\u adulescu proved, in
\cite{Ra2} Prop.2.1, that the norms $||\cdot||_R$ define a structure
of Frechet space on V, i.e. locally convex space, metrizable (with a
metric invariant by translations) and complete. In this case $V^*$
is separating for $V$ and thus we can consider the
$\sigma(V,V^*)$-topology on $V$.\\
Now we want to generalize the notion of ''square'' and ''sum of
squares''. Starting from the classical theory, in which the squares
are elements of the form $a^*a$, the first step is to define an
adjoint operation on $V$.
\begin{defin}
We set $(Y_{i_1}\cdot...\cdot Y_{i_p})^*=Y_{i_p}\cdot...\cdot Y_1$,
and $a^*=\overline a$ for the coefficients. We can extend this
mapping by linearity to an adjoint operation on $V$.
\end{defin}
Now, observing that series are too general to obtain in a finite
number of steps, we have quite naturally the following
\begin{defin}
We say that $q\in V$ is sum of squares if it is in the weak closure
of the set of the elements of the form $\sum p^*p$, $p\in V$.
\end{defin}
Now we observe that the formulation of Hilbert's 17th problem with
matrices regards matrices whose entries are REAL polynomial, the
geometric formulation regards analytic functions with REAL values.
So, recalling that REAL in operator algebras becomes SELF-ADJOINT,
we have that our natural setting to generalize Hilbert's 17th
problem is not $V$, but $V_{sa}=\{v\in V:v^*=v\}$.\\
It remains only to generalize the notion of positivity.
\begin{defin}
Let $p\in V_{sa}$. We say that $p$ is positive semidefinite if for
every $N\in\mathbb N$ and for every $X_1,...,X_n\in M_N(\mathbb C)$,
one has
$$
tr(p(X_1,...,X_n))\geq0
$$
$V_{sa}^+$ will denote the set of positive semidefinite elements of
$V_{sa}$.
\end{defin}
In order to arrive to the generalization of Hilbert's 17th problem
we have to do a last
\begin{rem}\label{noncommutative}
In the case of polynomials in $\mathbb R[Y_1,...,Y_n]$, we have
$Y_I-Y_{\tilde I}=0$, for any permutation $\tilde I$ of $I$. In our
non-commutative case, we cannot have this equality and thus we have
to identify elements which \emph{differ by permutation}. A way to do
this identification is given by the following
\begin{defin}
Two elements $p,q\in V_{sa}^+$ are called cyclic equivalent if $p-q$
is weak limit of sums of scalars multiples of monomials of the form
$Y_I-Y_{\tilde I}$, where $\tilde I$ is a cyclic permutation of $I$.
\end{defin}
\end{rem}
In this way, we have the following
\begin{prob}{\bf (Operator Algebra
version)} Is every element in $V_{sa}^+$ cyclic equivalent to a sum
of squares?
\end{prob}
Here is the beautiful and unexpected theorem by R\u adulescu.
\begin{teo}{\bf (R\u adulescu, \cite{Ra2} Cor.1.2)}
The following statements are equivalent
\begin{enumerate}
\item Connes' embedding conjecture is true.
\item Operator algebra version of Hilbert's 17th problem has
affirmative answer.
\end{enumerate}
\end{teo}
Following Radulescu some authors have began an approach to Connes'
embedding problem via \emph{sums of hermitian squares}. In this last
page we want to describe briefly the main result of
Klep-Schweighofer's work (see \cite{Kl-Sc} and also \cite{Ju-Po}
for a development).\\\\
Let $K$ be the real or the complex field and $V=K[Y_1,...,Y_n]$. So,
the first difference between this approach and Radulescu's one is
that Klep and Schweighofer work with polynomial and Radulescu works
with formal series. Other differences are given by the choice of the
adjoint operation and the cyclic equivalence. More precisely, they
take the identity operation (on the inderminates) as adjoint
operation and the following as equivalence
\begin{defin}
$p,q\in V$ are called equivalent if $p-q$ is sum of commutators. We
write $p\sim q$.
\end{defin}
Once again this equivalence relation is clearly trivial in the
commutative case.
\\\\
On the other hand, the notion of positivity introduced by Klep and
Schweighofer is a little less strong
\begin{defin}
$f\in V$ is called positive semidefinite if for any $s\in\mathbb N$
and for any contractions $A_1,...A_n\in M_s(\mathbb R)$ one has
$$
tr(f(A_1,...,A_n))\geq0
$$
The set of positive semidefinite element will be denoted by $V^+$.
\end{defin}
Now we give the definition of quadratic module, which is the major
difference with Radulescu's formulation.
\begin{defin}
A subset $M\subseteq V_{sa}$ is called quadratic module if the
followings hold
\begin{enumerate}
\item $1\in M$
\item $M+M\subseteq M$
\item $p^*Mp\subseteq M$, for all $p\in V$.
\end{enumerate}
The quadratic module generated in $V$ by the elements
$1-X_1^2,...,1-X_n^2$ will be denoted by $Q$.
\end{defin}
\begin{teo}{\bf (Klep-Schweighofer)}
The following statements are equivalent
\begin{enumerate}
\item Connes' embedding conjecture is true.
\item The following Radulescu's type implication holds
$$
f\in
V^+\,\,\,\,\,\,\,\,\Rightarrow\,\,\,\,\,\,\,\,\forall\varepsilon>0,
\exists q\in Q\,\,s.\,\,t.\,\,f+\varepsilon\sim q
$$
\end{enumerate}
\end{teo}

\subsection{Voiculescu's entropy}

In order to show the relation between Connes' embedding conjecture
and Voiculescu's free entropy, we have to recall briefly
Voiculescu's definition. References for this part are the
preliminary sections of the papers by Voiculescu \cite{Vo1} and
\cite{Vo2}. A motivation for these definitions can be found in
\cite{Vo3}.
\begin{nota}
We recall a construction of the entropy of a random variable which
outcomes the set $\{1,...n\}$ with probabilities $p_1,... p_n$. The
microstates are the set
$$
\{1,...n\}^N=\{f:\{1,...N\}\rightarrow\{1,...n\}\}
$$
The set of microstates which $\varepsilon$-approximate the discrete
distribution $p_1,...p_n$ is
$$
\Gamma(p_1,...
p_n,N,\varepsilon)=\{f\in\{1,...n\}^N:|\frac{|f^{-1}(i)|}{N}-p_i|<\varepsilon\,\,\,\,\forall
i=1,...n\}
$$
where $|f^{-1}(i)|$ is the number of elements in the counter-image.
Now, one takes the limit of
$$
N^{-1}lg|\Gamma(p_1,...p_n,N,\varepsilon)|
$$
as $N\rightarrow\infty$ and then lets $\varepsilon$ go to zero. Thus
we obtain the classical formula $\sum p_ilgp_i$ for the entropy.
\end{nota}
Voiculescu has generalized this construction to the non-commutative
setting of von Neumann algebras.
\begin{nt}{\bf (Lebesgue measure instead of the discrete one)}
Let $k$ be a positive integer and $(M_k(\mathbb C)_{sa})^n$ the set
of $n$-tuples of self-adjoint $k\times k$ complex matrices. Let
$\lambda$ be the Lebesgue measure on $(M_k(\mathbb C)_{sa})^n$
corresponding to the Euclidean norm
$$
||(A_1,...,A_n)||_{HS}^2=Tr(A_1^2+...+A_n^2)
$$
where $Tr$ is the non-normalized trace on $M_k(\mathbb C)$.
\end{nt}
\begin{nt}{\bf (Microstates are matrices)}
Fixed $\varepsilon, R>0$ and $m,k\in\mathbb N$. Let $X_1,...X_n$
free random variables on a finite factor $M$. We set
$$
\Gamma_R(X_1,...,X_n;m,k,\varepsilon)=\{(A_1,...,A_n)\in M_k(\mathbb
C)_{sa}^n\,\,\,s.t. $$ $$\left\{
                      \begin{array}{ll}
                        ||A_j||\leq R \\
                        |tr(A_{i_1}\cdot...\cdot
A_{i_p})-\tau(X_{i_1}\cdot...\cdot
X_{i_p})|<\varepsilon\,\,\,\forall (i_1,...,i_p)\in\{1,...n\}^p,
1\leq p\leq n
                      \end{array}
                    \right.
\}
$$
\end{nt}
\begin{defin}{\bf (Generalization of the process of limit)}
$$
\chi_R(X_1,...,X_n;m,k,\varepsilon)=log\lambda(\Gamma_R(X_1,...X_n;m,k,\varepsilon))
$$
$$
\chi_R(X_1,...X_n;m,\varepsilon)=limsup_{k\rightarrow\infty}(k^{-2}\chi_R(X_1,...X_n;m,k,\varepsilon)+2^{-1}nlog(k))
$$
$$
\chi_R(X_1,...X_n)=inf\{\chi_R(X_1,...X_n;m,\varepsilon),m\in\mathbb
N,\varepsilon>0\}
$$
$$
\chi(X_1,...X_n)=sup\{\chi_R(X_1,... X_n), R>0\}
$$
$\chi(X_1,... X_n)$ is called free entropy of the variables
$X_1,...X_n$.
\end{defin}
\begin{nota}
The factor $k^{-2}$ instead of $k^{-1}$ comes from the
normalization. The addend $2^{-1}nlg(k)$ is necessary, since
otherwise $\chi_R(X_1,... X_n;m,\varepsilon)$ should be always equal
to $-\infty$.
\end{nota}
By definition it follows that the free entropy can be equal to
$-\infty$. Voiculescu himself has found some examples
\begin{prop}{\bf (\cite{Vo2}, Prop. 3.6,c))}
If $X_1,... X_n$ are linearly dependent, then $\chi(X_1,...
X_n)=-\infty$.
\end{prop}
In order to have $\chi(X_1,...X_n)>-\infty$ we need at least that
$\Gamma_R(X_1,... X_n,m,k,\varepsilon)$ is not empty for some $k$,
i.e. the finite subset $X=\{X_1,...X_n\}$ of $M_{sa}$ has
microstates. This requirement is equivalent to Connes' embedding
conjecture:
\begin{teo}\label{microstates}
Let $M$ be a type $II_1$ factor. The following conditions are
equivalent
\begin{enumerate}
\item Every finite subsets $X\subseteq M_{sa}$ has microstates.
\item $M$ is embeddable into $R^\omega$.
\end{enumerate}
\begin{proof}
($1.\Rightarrow2.$)\\
Let $Y=\{x_1,x_2,...\}$ a norm-bounded generating set for $M$. Fix
$m\in\mathbb N$ and $\varepsilon=\frac{1}{m}$. By hypothesis, there
exist a natural number $k$ and $A_1^{(m)},...A_m^{(m)}\in
M_k(\mathbb C)$ which are microstates for $x_1,... x_m$. It has been
proved by Voiculescu (see \cite{Vo2}) that one can choose
$||A_i^{(m)}||\leq||x_i||$. Let $\pi_k:M_k(\mathbb C)\rightarrow R$
any unital *monomorphism. Define $a_i^m=\pi_k(A_i^{(m)})$ and
$b_i=\{a_i^m\}_{m=1}^\infty\in l^\infty(R)$, where $a_i^m=0$, if
$i>m$. Let $z_i$ be the image of $b_i$ into $R^\omega$. The mapping
$x_i\rightarrow z_i$ extends to an embedding $M\hookrightarrow
R^\omega$.\\\\
($2.\Rightarrow1.$)\\
Let $X=\{x_1,...x_n\}\subseteq R^\omega$. These elements are 2-norm
limit of element of $R$ and thus we can find $a_1,... a_n\in R$
whose mixed moments approximate those of the $x_i$'s (indeed
$\tau_{R^\omega}(x_{i_1}\cdot...\cdot
x_{i_p})=lim\tau_R(x_{i_1}^{(n)}\cdot...\cdot x_{i_p}^{(n)})$). Thus
the implication follows by noting that every finite subsets of $R$
is 2-norm approximately contained in some copy of $M_k(\mathbb R)$,
for some $k$ sufficiently large.
\end{proof}
\end{teo}

\subsection{Collins and Dykema's approach via eigenvalues}
Connes' embedding problem regards the approximation of the operators
in a separable type $II_1$ factor via matrices. The basic idea of
the approach by Collins and Dykema is that such an approximation
must reflect on the eigenvalues: the eigenvalues of an operator in a
separable type $II_1$ factor should be approximated by the
eigenvalues of the matrices. This is just the basic idea, but there
are some problems:
\begin{enumerate}
\item What does \emph{eigenvalue} mean for an operator in a
separable type $II_1$ factor?
\item In which sense those \emph{eigenvalues} are approximated by
the eigenvalues of the matrices?
\end{enumerate}
We start by answering to the first question.\\
Let $M$ be a separable type $II_1$ factor and $\tau$ its unique
faithful normalized trace. For any $a\in M_{sa}$ we can define the
distribution of $a$ as the Borel measure $\mu_a$, supported on the
spectrum of $a$, such that
$$
\tau(a^n)=\int_{\mathbb
R}t^nd\mu_a(t)\,\,\,\,\,\,\,\,\,\,\,\,\,\,\,\,\,\, n\geq1
$$
\begin{defin}
Let $a\in M_{sa}$. The eigenvalue function of $a$ is the function
$\lambda_a:[0,1)\rightarrow\mathbb R$ defined by
$$
\lambda_a(t)=sup\{x\in\mathbb R:\mu_a((x,+\infty))>t\}
$$
\end{defin}
This definition generalizes what happens in $M_N(\mathbb C)$ as
follows:\\\\
Let $a\in M_N(\mathbb C)_{sa}$ and let
$\alpha=(\alpha_1,...\alpha_N)$ be its eigenvalue sequence, i.e.
$\alpha_1,....\alpha_N$ are the eigenvalues of $a$ listed in
non-increasing order and according to their multiplicity. In this
case one has
$$
\lambda_a(t)=\alpha_j
$$
where $j$ is defined by the property $\frac{j-1}{N}\leq
t<\frac{j}{N}$.\\\\
Now we pass to the second question. First of all we need a topology
with respect to we can consider the approximations. We denote
$$
\mathcal F=\{f:[0,1)\rightarrow\mathbb
R\,\,\,\,right-continuous,\,\,\,non-increasing\,\,\,and\,\,\,bounded\}
$$
Clearly any eigenvalue function belongs into $\mathcal F$.
Conversely, given $f\in\mathcal F$ and $M\in\Im_{II_1}$, there
exists $a\in M_{sa}$ such that $\lambda_a=f$. In this way we are
able to identify $\mathcal F$ with the set of eigenvalue functions.
On the other hand
$$
\mathcal
F=\{compactly-supported\,\,\,Borel\,\,\,measures\,\,\,on\,\,\,\mathbb
R\}\subseteq C(\mathbb R)^*
$$
Therefore a natural topology on $\mathcal F$ (and thus on the set of
eigenvalue functions too) is the weak* topology on $C(\mathbb
R)^*$.\\\\
Above we said that the notion of eigenvalue function for operators
generalizes that for matrices. Now we need to give a little
formalization of this fact. Let $\mathbb R_\leq^N$ be the set of
$N$-tuples of real numbers listed in non-increasing order. The
correspondence
$$
\alpha=(\alpha_1,...\alpha_N)\in\mathbb
R_\leq^N\rightarrow\lambda_\alpha(t)=\alpha_j\,\,\,\,\,\,where\,\,\,\,\,\frac{j-1}{N}\leq
t<\frac{j}{N}
$$
gives an embedding $\mathbb R_\leq^N\subseteq\mathcal F$. This
embedding is very good, since it preserves the affine structure (the
affine structure on $\mathcal F$ is defined by taking the usual
scalar multiplication and sum of functions; the affine structure on
$\mathbb R_\leq^N$ comes from $\mathbb R^N$).\\\\
Now the idea is that Connes' embedding conjecture should be
equivalent in something like the density of $\mathbb R_\leq^N$ into
$\mathcal F$. Actually it happens something more precise and
elegant. In two words: Connes' embedding conjecture is equivalent to
the possibility of approximating the eigenvalue function of
operators of the form
$$
a_1\otimes x_1+a_2\otimes x_2\,\,\,\,\,\,\,\,\,\,\,\,\,\, a_i\in
M_N(\mathbb C), x_i\in M\,\,\,\,\,\,\,\, (M\in\Im_{II_1})
$$
with the eigenvalue function of operators of the form
$$
a_1\otimes y_1+a_2\otimes y_2\,\,\,\,\,\,\,\,\,\,\,\,\,\,\
a_i,y_i\in M_N(\mathbb C)
$$
where the eigenvalues functions of the $y_i$'s are the same of those
of the $x_i$'s (after the embedding $\mathbb
R_\leq^N\subseteq\mathcal F)$.\\\\
We give some details in order to arrive to the correct enunciation
of Collins-Dykema's theorem. Let $\alpha,\beta\in\mathbb R_\leq^N,
d\in\mathbb N,a_1,a_2\in M_{Nd}(\mathbb C)_{sa}, M\in\Im_{II_1}$. We
denote
$$
K_{\alpha,\beta,d}^{a_1,a_2}=\{\lambda_C, C=a_1\otimes
U(diag(\alpha)\otimes Id_d)U^*+a_2\otimes V(diag(\beta)\otimes
Id_d)V^*, U,V\in U(M_{nd}(\mathbb C))\}
$$
$$
K_{\alpha,\beta,\infty}^{a_1,a_2}=\overline{\bigcup_{d\in\mathbb
N}K_{\alpha,\beta,d}^{a_1,a_2}}
$$
where the closure is respect to the weak* topology on $\mathcal F$.
$$
L_{\alpha,\beta,M}^{a_1,a_2}=\{\lambda_C, C=a_1\otimes
x_1+a_2\otimes x_2\}
$$
where $x_1,x_2\in M$ whose eigenvalue functions agree with those of
the matrices $diag(\alpha)$ and $diag(\beta)$.\\ At last we denote
$$
L_{\alpha,\beta}^{a_1,a_2}=\bigcup_{M\in\Im_{II_1}}L_{\alpha,\beta,M}^{a_1,a_2}
$$
Here is Collins-Dykema's theorem
\begin{teo}{\bf (Collins-Dykema, \cite{Co-Dy}, Th. 4.6)}
The following statements are equivalent
\begin{enumerate}
\item Connes' embedding conjecture is true.
\item $L_{\alpha,\beta}^{a_1,a_2}=K_{\alpha,\beta,\infty}^{a_1,a_2}$
\end{enumerate}
\begin{proof}
\emph{(Sketch)}. If Connes' embedding conjecture is true, then
$L_{\alpha,\beta}^{a_1,a_2}=L_{\alpha,\beta, R^\omega}^{a_1,a_2}$.
On the other hand $L_{\alpha,\beta,
R^\omega}^{a_1,a_2}=K_{\alpha,\beta,\infty}^{a_1,a_2}$. Hence the
first implication easily follows. Conversely, one can suppose that
$M$ is generated by two self-adjoint elements $x_1,x_2$.
Approximating $x_1,x_2$ we can assume their eigenvalue function
belong into $\mathbb R_\leq^N$, for some $N$. By adding constants we
may also assume that $x_1,x_2$ are positive and invertible. Now a
theorem by Collins and Dykema (see \cite{Co-Dy}, 3.6) shows that
$x_1,x_2$ have microstates and thus the thesis follows from
Th.\ref{microstates}.
\end{proof}
\end{teo}
\section{Acknowledgements}
I want to thank my supervisor Florin R\u adulescu and my colleague
and friend Liviu Paunescu for their useful suggestions and
corrections, Profs. N.P. Brown, B. Collins and K. Dykema who have
helped me to develop the section regarding their paper.

VALERIO CAPRARO, \emph{Universit\`{a} degli Studi di Roma "Tor
Vergata"}, email: capraro@mat.uniroma2.it\\\\
All comments are very much welcome! I will be happy to correct
inaccuracies or omissions in future versions.

\end{document}